\documentclass[hidelinks,onefignum,onetabnum]{siamart251216}
\usepackage{amsfonts}
\usepackage{amssymb}
\usepackage{bm}
\usepackage{booktabs}
\usepackage{color}
\definecolor{RED}{rgb}{1,0,0}
\definecolor{ORANGE}{rgb}{1,0.45,0}

\headers{Strong NP-Completeness and Exact Algorithms}
{Linbin Li, Yong Xia, and Zhuoyi Xu}
\title{Maximin Dispersion over Simple Polyhedral Domains:\\
Strong NP-Completeness and Exact Algorithms\footnotemark[1]\thanks{Submitted on Nov. 16th, 2026.
\funding{This work was supported by the National Natural Science Foundation of China (No.
12631012 and No. 12401410).}}}
\author{Linbin Li\footnotemark[2]\thanks{School of Mathematical Sciences, Beihang University,
Beijing, China} (\email{linbinli@buaa.edu.cn, yxia@buaa.edu.cn}).
\and Yong Xia\footnotemark[2]
\and Zhuoyi Xu\footnotemark[3]\thanks{School of Statistics and Data Science, Capital University of Economics and Business,
Beijing, China} (\email{xuzhuoyi@cueb.edu.cn}).}

\ifpdf
\hypersetup{
  pdftitle={Maximin Dispersion over the Unit $\ell_1$-Ball and Vertex-Generated Simplices: Strong NP-Completeness and Exact Algorithms},
  pdfauthor={Linbin Li, Yong Xia, and Zhuoyi Xu}
}
\fi

\begin{document}

\maketitle

\begin{abstract}
In this paper, we prove that the decision version of weighted maximin dispersion over the unit $\ell_1$-ball
with rational sites and positive rational weights is strongly NP-complete.
The same holds for unweighted dispersion over full-dimensional integer
simplices when the sites are exactly their vertices.
Both reductions from Exact Cover by 3-Sets rely on a shared quadratic
rigidity argument.
We also give exact algorithms that run in polynomial time for every fixed
dimension.
For the simplex problem, a Voronoi decomposition yields rational optimal
solutions.
For the weighted $\ell_1$-ball problem, an orthant decomposition and a
single-paraboloid lifting yield optimal solutions whose coordinates lie in a common real algebraic extension of degree at most two.
In both cases, these solutions have polynomial encoding length and provide
certificates for NP membership.
\end{abstract}



\section{Introduction}
Given a compact set $P\subseteq\mathbb R^d$, sites
$\mathcal A=\{\bm a_1,\ldots,\bm a_r\}\subseteq\mathbb R^d$, and positive
weights $\bm\omega=(\omega_1,\omega_2,\ldots,\omega_r)^\top\in\mathbb R^r_{+}$, weighted maximin dispersion asks for
\begin{equation}
 \operatorname{MD}(P,\mathcal A,\bm\omega)
 =\max_{\bm x\in P}\min_{1\leq s\leq r}
       \omega_s\|\bm x-\bm a_s\|_2^2.
 \label{eq:general-wmd}
\end{equation}
It models the placement of a point as far as possible from prescribed sites
in a weighted worst-case sense.  Related maxmin location, distance-design,
and dispersion models arise in facility location, statistics, and
multiobjective selection
\cite{DasarathyWhite1980,JohnsonEtAl1990,RaviEtAl1994,White1996}. Convex relaxations and exactness conditions for weighted dispersion are
developed in \cite{HainesEtAl2013,WangXia2016,JeyakumarLi2018};
approximation methods include random sampling
\cite{WangXu2021,WuXiaWang2018} and first-order optimization
\cite{XuWangHuang2021}.  Hardness over elementary domains also occurs in
other nonconvex optimization models
\cite{PardalosVavasis1991,deKlerk2008}.

One might expect the complexity of \cref{eq:general-wmd} to be governed
primarily by the geometry of $P$.  Existing results show, however, that
hardness persists even over elementary domains.  Haines et al.\
\cite{HainesEtAl2013} established NP-hardness for the box-constrained
problem, including the unit $\ell_\infty$-ball, and Wang and Xia
\cite{WangXia2016} proved NP-hardness for the unit $\ell_2$-ball; both
results allow equal weights.  In these models the norm specifies the
feasible region, while the objective continues to use squared Euclidean
distances.  This leaves the third canonical norm ball, the unit
$\ell_1$-ball, as a natural missing case.

We first consider the specialization of \cref{eq:general-wmd} in which the
feasible region is fixed, for each input dimension $d$, to the unit
$\ell_1$-ball
\begin{equation}
 B_1^d
 :=
 \left\{
   \bm{x}\in\mathbb{R}^d:
   \|\bm{x}\|_1\leq 1
 \right\}.
 \label{eq:intro-l1-ball}
\end{equation}
Given rational sites
$\mathcal{A}=\{\bm{a}_1,\ldots,\bm{a}_r\}\subseteq\mathbb{Q}^d$
and positive rational weights
$\bm{\omega}=(\omega_1,\ldots,\omega_r)^{\mathsf T}$, define weighted maximin dispersion over the unit $\ell_1$-ball by
\begin{equation}
 (\ell_1\text{-WMD}) \quad \Phi^*(\mathcal{A},\bm{\omega})
 :=\max_{\bm{x}\in B_1^d} F(\bm{x}):=
 \max_{\bm{x}\in B_1^d}
 \min_{1\leq s\leq r}
 \omega_s\|\bm{x}-\bm{a}_s\|_2^2.
 \label{eq:intro-l1-wmd}
\end{equation}
Given a rational threshold $\tau\geq0$, we refer to the decision version of weighted maximin dispersion over the
unit $\ell_1$-ball as $\ell_1$-WMD-DEC.
\begin{equation}
 (\ell_1\text{-WMD-}\text{DEC})\quad\Phi^*(\mathcal{A},\bm{\omega})\geq\tau.
 \label{eq:intro-l1-decision}
\end{equation}

We prove that $\ell_1$-WMD-DEC is strongly NP-complete with positive
rational weights.  We then complement the hardness result with an exact
algorithm: after decomposing the $\ell_1$-ball into orthant simplices, a
single lifting variable makes every weighted distance constraint linear and
leaves only one common paraboloid equation.  This structure also yields a
polynomial-length quadratic algebraic certificate.

The nonnegative orthant of the unit $\ell_1$-ball is a simplex, suggesting a
second, structurally different restriction.  We ask whether hardness persists
when the sites are no longer independent data but must coincide with the
vertices generating the feasible simplex.  In the $\ell_1$-ball model the
domain is fixed and only the sites and weights vary; in the simplex model the
domain and the sites are coupled.  The two restrictions therefore require
different geometric realizations of a shared exact-cover rigidity mechanism.

For the second model, let
\begin{equation}
 S
 =
 \operatorname{conv}
 \{\bm{v}_0,\ldots,\bm{v}_m\}
 \subseteq\mathbb{R}^d,
 \qquad m\leq d,
 \label{eq:intro-svbc-simplex}
\end{equation}
be a nondegenerate rational $m$-simplex.  The simplex vertex ball
coverage problem (SVBC) computes
\begin{equation}
 (\text{SVBC})\quad R^*(S)
 :=
 \max_{\bm{x}\in S}
 \min_{0\leq k\leq m}
 \|\bm{x}-\bm{v}_k\|_2^2.
 \label{eq:intro-svbc}
\end{equation}
Equivalently, $R^*(S)$ is the smallest common squared radius of
closed balls centered at the vertices that cover $S$.
For a rational threshold $\tau\geq0$, we define SVBC-DEC to ask whether
\begin{equation}
 (\text{SVBC-DEC})\quad R^*(S)\geq\tau.
 \label{eq:intro-svbc-decision}
\end{equation}
We use the dispersion decision convention; coverage by closed balls
of squared radius $\tau$ instead corresponds to $R^*(S)\leq\tau$.

Vertex-centered simplex covering was studied by Casado et al.\
\cite{CasadoEtAl2011} using power diagrams; for a common radius,
the critical squared covering radius is $R^*(S)$.
Zhang and Xia \cite{ZhangXia2022} proved NP-hardness for simplex
covering optimization with independently specified centers and
radii, and gave a strongly polynomial algorithm for a fixed
number of balls.
Here we establish strong NP-completeness under the restriction
that the sites are exactly the vertices of a full-dimensional integer simplex. Although the fixed-ball algorithm already implies polynomial
solvability of fixed-dimensional SVBC, we give a direct
Voronoi-based enumeration with an explicit candidate bound
and polynomial-length rational certificates for NP membership.

\subsection{Main results}
The two principal complexity results are the following.

\begin{theorem}
\label{thm:l1-main}
The problem $\ell_1$-WMD-DEC is strongly NP-complete.
\end{theorem}

\begin{theorem}
\label{thm:main}
The problem SVBC-DEC is strongly NP-complete even for full-dimensional integer simplices.
\end{theorem}

The paper first develops the two hardness reductions and then
their complementary exact algorithms.
\Cref{sec:preliminaries} introduces the shared exact-cover rigidity
lemma.  \Cref{sec:l1-wmd} realizes its conditions through weighted
sites over the unit $\ell_1$-ball, and \cref{sec:svbc-hardness}
realizes them through the vertices of a full-dimensional integer
simplex.  \Cref{sec:l1-exact,sec:voronoi} give the respective exact
algorithms and establish NP membership, completing
\cref{thm:l1-main,thm:main}.
The auxiliary quadratic computations and the rational cluster
construction are proved in
\cref{app:one-quadric-details,app:svbc-realization}.

Throughout, bold lowercase and uppercase letters denote vectors and
matrices, respectively.  For a positive integer $k$, let
$\bm{0}_k,\bm{1}_k\in\mathbb{R}^k$ denote the zero and all-ones vectors,
and let $\bm{e}_i$ denote the $i$th standard basis
vector.  
All rational input data are encoded in binary.

\subsection{A shared exact-cover rigidity lemma }
\label{sec:preliminaries}

Both hardness reductions are from Exact Cover by 3-Sets
(X3C). An instance of X3C consists of a universe
$U=\{1,\ldots,3q\}$ and a family
$\mathcal{S}=\{S_1,\ldots,S_n\}$ with $|S_i|=3$ for every $i$.
The problem asks whether some $q$ sets in $\mathcal{S}$ partition $U$.
X3C is strongly NP-complete \cite[Problem SP2]{GareyJohnson1979}.

Let
\begin{equation}
 \bm{B}\in\{0,1\}^{3q\times n},
 \qquad
 B_{ji}=1
 \quad\Longleftrightarrow\quad
 j\in S_i,
 \label{eq:incidence-matrix}
\end{equation}
and let $\bm{b}_i\in\{0,1\}^{3q}$ denote the $i$th column of
$\bm{B}$.  Since every set contains three elements,
\begin{equation}
 \bm{B}^{\mathsf T}\bm{1}_{3q}=3\bm{1}_n.
 \label{eq:column-sums}
\end{equation}
The following lemma recovers an exact cover from two quadratic equalities on the simplex and supplies the common rigidity argument for both reductions.
\begin{lemma}[Exact-cover rigidity]
\label{lem:exact-cover-rigidity}
Suppose that $\bm z\in\mathbb R_+^n$ satisfies
\begin{equation}
 \bm{1}_n^{\mathsf T}\bm{z}=1,
 \qquad
 \|\bm{z}\|_2^2=\frac{1}{q},
 \qquad
 \|\bm{B}\bm{z}\|_2^2=\frac{3}{q}.
 \label{eq:rigidity-hypotheses}
\end{equation}
Then $ q\bm{z}\in\{0,1\}^n$ and $ \bm{B}(q\bm{z})=\bm{1}_{3q}$.
Consequently, the support of $\bm{z}$ indexes an exact cover.
\end{lemma}

\begin{proof}
By \cref{eq:column-sums} and the Cauchy--Schwarz inequality, $\|\bm B\bm z\|_2^2
 \geq\frac{(\bm1_{3q}^{\mathsf T}\bm B\bm z)^2}{3q}
 =\frac3q $.
Equality in \cref{eq:rigidity-hypotheses} therefore forces
\begin{equation}
 \bm B\bm z=\frac1q\bm1_{3q}.
 \label{eq:shared-uniform-incidence}
\end{equation}
For every $i$, choose any $j\in S_i$.  Nonnegativity and
\cref{eq:shared-uniform-incidence} give $0\leq z_i\leq1/q$.  Hence
\begin{equation}
 0=\frac1q-\|\bm z\|_2^2
  =\sum_{i=1}^n z_i\left(\frac1q-z_i\right).
 \label{eq:shared-coordinate-rigidity}
\end{equation}
Every summand is nonnegative, so every $z_i$ belongs to
$\{0,1/q\}$.  The coordinates sum to one, and
\cref{eq:shared-uniform-incidence} proves the exact-cover claim.
\end{proof}

Throughout the reductions, a YES-instance of X3C admits an exact cover, whereas a NO-instance does not. Both reductions enforce the hypotheses of
\cref{lem:exact-cover-rigidity} through distance constraints.
In the $\ell_1$-ball reduction, $\bm z$ lies on the nonnegative
simplex face; in the SVBC reduction, it consists of barycentric
coefficients on the candidate face.
Each construction gives an upper bound $\tau$ on the
dispersion objective.  An exact cover attains $\tau$, whereas a point
attaining $\tau$ yields an exact cover by
\cref{lem:exact-cover-rigidity}.  Thus the optimal value equals $\tau$
for YES-instances and is strictly smaller for NO-instances.
To establish strong NP-hardness, we also bound the number of output
entries and the numerical magnitudes of all rational numerators and
denominators by polynomials in the explicit X3C input size.

\section{Strong NP-hardness over the unit
\texorpdfstring{$\ell_1$}{l1}-ball}
\label{sec:l1-wmd}

This section proves the strong NP-hardness part of \cref{thm:l1-main}.
The first step shows that, when all sites lie in the strict negative orthant, the unit $\ell_1$-ball can
be reduced to its probability-simplex face. We then construct a weighted
site family whose distance inequalities force the hypotheses of
\cref{lem:exact-cover-rigidity}.

\subsection{Reduction to the nonnegative simplex face}

\begin{lemma}
\label{lem:l1-folding}
Suppose that every site lies in the strict negative orthant; that is,
$a_{sj}<0$ for all $1\leq s\leq r$ and $1\leq j\leq d$.
Then every maximizer of \cref{eq:intro-l1-wmd} belongs to
\begin{equation}
 \Delta_{d-1}^{+}
 :=
 \left\{
   \bm{x}\in\mathbb{R}^d:
   \bm{x}\geq\bm{0}_d,\ 
   \bm{1}_d^{\mathsf T}\bm{x}=1
 \right\}.
 \label{eq:l1-positive-face}
\end{equation}
Consequently,
\begin{equation}
 \max_{\|\bm{x}\|_1\leq1}
 \min_{1\leq s\leq r}
 \omega_s\|\bm{x}-\bm{a}_s\|_2^2
 =
 \max_{\bm{x}\in\Delta_{d-1}^{+}}
 \min_{1\leq s\leq r}
 \omega_s\|\bm{x}-\bm{a}_s\|_2^2.
 \label{eq:l1-folded-value}
\end{equation}
\end{lemma}

\begin{proof}
We move any point outside the nonnegative simplex face to that
face while increasing its distance from every site.
Fix $\bm x\in B_1^d$ and let
\[
 \bm u:=|\bm x|,\qquad
 \varepsilon:=1-\|\bm x\|_1,\qquad
 \widehat{\bm x}:=\bm u+\varepsilon\bm e_k^{(d)},
\]
where $|\bm x|$ is taken componentwise and
$k\in\{1,\ldots,d\}$ is arbitrary.
Then $\varepsilon\geq0$ and
$\widehat{\bm x}\in\Delta_{d-1}^{+}$.
For every site $\bm a_s$,
\begin{align}
 \|\bm u-\bm a_s\|_2^2-\|\bm x-\bm a_s\|_2^2
 &=4\sum_{j:x_j<0}|x_j|\,|a_{sj}|,
 \label{eq:l1-folding-identity}\\
 \|\widehat{\bm x}-\bm a_s\|_2^2-\|\bm u-\bm a_s\|_2^2
 &=2\varepsilon(u_k-a_{sk})+\varepsilon^2.
 \label{eq:l1-boundary-increase}
\end{align}
Both differences are nonnegative.  If
$\bm x\notin\Delta_{d-1}^{+}$, then either $\bm x$ has a
negative component or $\varepsilon>0$, so at least one
difference is strictly positive for every $s$.
Since the site family is finite and all weights are positive,
$F(\widehat{\bm x})>F(\bm x)$ whenever
$\bm x\notin\Delta_{d-1}^{+}$, where
\(
 F(\bm x):=\min_{1\leq s\leq r}
 \omega_s\|\bm x-\bm a_s\|_2^2
\)
is the objective function in \cref{eq:intro-l1-wmd}.
Thus every maximizer lies in $\Delta_{d-1}^{+}$, proving
\cref{eq:l1-folded-value}.
The same transformation preserves threshold feasibility,
so a threshold-feasible point exists in $B_1^d$ if and only
if one exists in $\Delta_{d-1}^{+}$.
\end{proof}

\subsection{Reduction from X3C}

\begin{proof}[Proof of the strong NP-hardness part of
\cref{thm:l1-main}]
Let
\[
 U=\{1,\ldots,3q\},
 \quad
 \mathcal{S}=\{S_1,\ldots,S_n\}, \quad \bm{B}\in\{0,1\}^{3q\times n}
\]
be an X3C instance and its incidence matrix, as defined in \cref{sec:preliminaries}. $B_{ji}=1$ if and only if $ j\in S_i$.
For each $j\in U$, let
$\bm b^{(j)}=(B_{j1},\ldots,B_{jn})^{\mathsf T}$ be the transpose
of the $j$th row of $\bm B$, and let
$r_j=\bm1_n^{\mathsf T}\bm b^{(j)}$ count the sets containing $j$. 

\textbf{Construction.}
The constructed WMD instance has dimension $d=n$ and threshold
$\tau=1$.  The auxiliary site $\bm a^0$ and the candidate sites
$\bm a^i$ will enforce $\|\bm x\|_2^2=1/q$, while the two sites
$\bm a^{j,+}$ and $\bm a^{j,-}$ for each element $j$ will enforce
$\bm B\bm x=\bm1_{3q}/q$.  We choose reciprocal weights to express
all these requirements at the same threshold.
Define the sites, thresholds, and weights by
\begin{align}
 \bm a^0&=-2\bm1_n,
 & T_0&=4n+4+1/q,
 & \omega_0&=1/T_0,
 \label{eq:l1-energy-site}\\
 \bm a^i&=\bm e_i^{(n)}-2\bm1_n,
 & T_c&=4n+1-1/q,
 & \omega_c&=1/T_c,
 \label{eq:l1-candidate-sites}\\
 \bm a^{j,+}&=\bm b^{(j)}-2\bm1_n,
 & T_j^+&=4n-3r_j+4-1/q,
 & \omega_j^+&=1/T_j^+,
 \label{eq:l1-plus-sites}\\
 \bm a^{j,-}&=-\bm b^{(j)}-2\bm1_n,
 & T_j^-&=4n+5r_j+4+3/q,
 & \omega_j^-&=1/T_j^-,
 \label{eq:l1-minus-sites}
\end{align}
where $1\leq i\leq n$ and $1\leq j\leq3q$; all candidate sites
$\bm a^i$ have weight $\omega_c$.
All denominators $T_0,T_c,T_j^+,T_j^-$ are positive.  Every site coordinate belongs to $\{-3,-2,-1\}$.  Hence all constructed sites lie in the
strict negative orthant, and \cref{lem:l1-folding} applies.

\textbf{Distance identities.}
We record how the constructed distances depend on the squared norm
and the incidence coordinates; these identities will verify both
directions of the reduction.
Fix $\bm x\in\Delta_{n-1}^{+}$ and set
$A=\|\bm x\|_2^2$ and $y_j=(\bm b^{(j)})^{\mathsf T}\bm x$.
Since $\bm1_n^{\mathsf T}\bm x=1$, direct expansion gives
\begin{align}
 D_0(\bm x)&:=\|\bm x-\bm a^0\|_2^2=A+4n+4,
 \label{eq:l1-distance-energy}\\
 D_i(\bm x)&:=\|\bm x-\bm a^i\|_2^2=A+4n+1-2x_i,
 \label{eq:l1-distance-candidate}\\
 D_j^+(\bm x)&:=\|\bm x-\bm a^{j,+}\|_2^2
              =A+4n-3r_j+4-2y_j,
 \label{eq:l1-distance-plus}\\
 D_j^-(\bm x)&:=\|\bm x-\bm a^{j,-}\|_2^2
              =A+4n+5r_j+4+2y_j,
 \label{eq:l1-distance-minus}
\end{align}
where $1\leq i\leq n$ and $1\leq j\leq3q$.

Each weight is the reciprocal of its corresponding threshold
($T_0$, $T_c$, $T_j^+$, or $T_j^-$).  At threshold \(\tau=1\), \(\bm x\) is feasible if and only if the following distance inequalities hold:
\begin{align}
 D_0(\bm{x})\geq T_0, \quad 
 D_i(\bm{x})\geq T_c, \quad
 D_j^+(\bm{x})\geq T_j^+,
 \quad
 D_j^-(\bm{x})\geq T_j^-.
\end{align}

\textbf{Completeness.}
An exact cover supplies a point with precisely the required norm
and incidence coordinates.
Suppose that $C\subseteq\{1,\ldots,n\}$ indexes an exact cover.
Define its normalized indicator vector $\bm\lambda^C\in\mathbb R^n$ by
\begin{equation}
 \lambda_i^C
 :=
 \begin{cases}
  1/q, & i\in C,\\
  0,   & i\notin C.
 \end{cases}
 \label{eq:l1-cover-point}
\end{equation}
Since $C$ consists of $q$ sets and covers every element exactly
once, $\bm\lambda^C\in\Delta_{n-1}^{+}$ and
\[
 A=\|\bm\lambda^C\|_2^2=\frac1q,
 \qquad
 y_j=\frac1q
 \quad(1\leq j\leq3q).
\]
Substituting into
\cref{eq:l1-distance-energy,eq:l1-distance-candidate,eq:l1-distance-plus,eq:l1-distance-minus} gives
\[
 D_0(\bm\lambda^C)=T_0,\qquad
 D_i(\bm\lambda^C)=T_c+\frac2q\mathbf1_{\{i\notin C\}},
 \qquad
 D_j^\pm(\bm\lambda^C)=T_j^\pm.
\]
Thus every weighted squared distance is at least one, and
the constructed instance is a YES-instance.

\textbf{Soundness.}
We now recover an exact cover from any threshold-feasible point.
The constraints associated with $\bm a^0$ and $\bm a^i$ first
determine $\|\bm x\|_2^2$; those associated with $\bm a^{j,+}$ and
$\bm a^{j,-}$ then determine $\bm B\bm x$.
Suppose that the constructed instance is a YES-instance.
By \cref{lem:l1-folding}, there is a point
$\bm x\in\Delta_{n-1}^{+}$ for which every weighted squared distance
is at least one.
The inequality $\omega_0\|\bm x-\bm a^0\|_2^2\ge 1$ and
\cref{eq:l1-distance-energy} imply
\begin{align}
 D_0(\bm{x})\geq T_0
 &\quad\Longrightarrow\quad
 A\geq\frac{1}{q}.
 \label{eq:l1-energy-lower}
\end{align}
Since $D_i(\bm x)\geq T_c$, $x_i\geq0$, and $\sum_i x_i=1$,
averaging the candidate constraints with weights $x_i$ gives
\begin{equation}
 T_c
 \leq\sum_{i=1}^n x_iD_i(\bm x)
 =\sum_{i=1}^n x_i(A+4n+1-2x_i)
 =4n+1-A.
 \label{eq:l1-candidate-average-identity}
\end{equation}
Using $T_c=4n+1-1/q$, we obtain $A\leq1/q$.
Together with \cref{eq:l1-energy-lower}, this yields $A=1/q$.
Substituting $A=1/q$ into
\cref{eq:l1-distance-plus,eq:l1-distance-minus} gives,
for every $j$,
\begin{align}
    0\leq
 D_j^+(\bm{x})-T_j^+ \notag =
 2\left(\frac{1}{q}-y_j\right), \quad 0
 \leq
 D_j^-(\bm{x})-T_j^- =2\left(y_j-\frac{1}{q}\right).
\end{align}
Thus for $1\leq j\leq3q $, one has $y_j={1}/{q}$.
Equivalently,
\begin{equation}
 \bm{B}\bm{x}
 =
 \frac{1}{q}\bm{1}_{3q},
 \qquad
 \|\bm{B}\bm{x}\|_2^2
 =
 3q\left(\frac{1}{q}\right)^2
 =
 \frac{3}{q}.
 \label{eq:l1-incidence-equality}
\end{equation}
Thus $\bm x$ satisfies the hypotheses of
\cref{lem:exact-cover-rigidity}, so its support indexes an exact cover.

We have proved the threshold equivalence
\begin{equation}
 \Phi^*(\mathcal{A},\bm{\omega})\geq1
 \quad\Longleftrightarrow\quad
 (U,\mathcal{S})\text{ has an exact cover}.
 \label{eq:l1-threshold-equivalence}
\end{equation}

The objective is also bounded above by one.  On
$\Delta_{n-1}^{+}$, if $A\leq1/q$, then $D_0(\bm x)\leq T_0$;
if $A>1/q$, then
$\sum_i x_iD_i(\bm x)=4n+1-A<T_c$, so some $D_i(\bm x)<T_c$.
The folding lemma extends this bound to the whole ball.

\textbf{Polynomial construction and strong hardness.}
The reduction produces one auxiliary site $\bm a^0$, $n$ candidate
sites $\bm a^i$, and two element sites $\bm a^{j,+},\bm a^{j,-}$
for each of the $3q$ elements.  Thus the indexed site
family contains $n+6q+1 $ sites in dimension $n$.

The weights can be written as
\begin{align}
 \omega_0
 &=
 \frac{q}{q(4n+4)+1},
 &
 \omega_c
 &=
 \frac{q}{q(4n+1)-1},
 \label{eq:l1-encoded-energy-candidate}
 \\
 \omega_j^+
 &=
 \frac{q}{q(4n-3r_j+4)-1},
 &
 \omega_j^-
 &=
 \frac{q}{q(4n+5r_j+4)+3}.
 \label{eq:l1-encoded-element-weights}
\end{align}
All numerators and denominators are positive.  Since
$0\leq r_j\leq n$, their magnitudes are at most $ q(9n+4)+3$, up to reduction of the fractions.  The threshold equals one, and every
site coordinate belongs to $\{-3,-2,-1\}$.  The number of sites, the
dimension, and the numerical magnitudes of all rational input data are
therefore polynomially bounded in the X3C input size.

Consequently, the output remains polynomial in length even if all
rational numerators and denominators are written in unary.  Since X3C is
strongly NP-complete, \cref{eq:l1-threshold-equivalence} proves that
$\ell_1$-WMD-DEC is strongly NP-hard.
\end{proof}

\section{Strong NP-hardness on a vertex-generated simplex}
\label{sec:svbc-hardness}
We next prove the strong NP-hardness part of \cref{thm:main}.
Building on the exact-cover encoding of \cref{sec:l1-wmd},
we introduce additional dimensions for the geometric construction.
Here the sites are unweighted and must coincide exactly with
the vertices of the feasible simplex.

The proof proceeds in three stages: a reduction on an embedded
simplex, a full-dimensional extension preserving threshold
equivalence, and a rational realization permitting integer scaling.

\subsection{Hardness on an embedded simplex}\label{sec:embedded-reduction}
Fix the X3C instance from \cref{sec:preliminaries}.
We may assume $q\geq3$.  When $q<3$, the X3C instance can be decided
directly in polynomial time and mapped to a fixed YES- or NO-instance
of SVBC-DEC.
Recall that \(n=|\mathcal S|\) and \(|U|=3q\).  The simplex constructed
in this subsection has intrinsic dimension \(m=2n\) and is embedded in
the ambient space \(\mathbb R^{d_0}\), where \(d_0=5n+3q+4\). 

\textbf{Construction.}
The vertices have three roles: candidate vertices $\bm p_i$
encode the X3C sets, a global anchor $\bm a$ helps confine
threshold-feasible points to the candidate face, and paired
anchors $\bm w_i$ separate the two quadratic quantities required
by \cref{lem:exact-cover-rigidity}.
Orthogonal coordinate blocks let us control these effects
through explicit squared-distance identities.
Set
\begin{equation}
 \bm u=\frac1q\bm1_{3q},\qquad
 L^2=1-\frac2q,\qquad
 H^2=4-\frac5q,\qquad
 \rho_q=4-\frac4q,\qquad
 K=H^2-\frac3q.
 \label{eq:svbc-parameters}
\end{equation}
Choose rational vectors $\bm g,\bm h\in\mathbb Q^4$ satisfying
$\|\bm g\|_2^2=L^2$ and $\|\bm h\|_2^2=H^2$.
Their polynomial-time construction is included in
\cref{lem:svbc-realization}, proved in
\cref{app:svbc-realization}.   For $1\leq i\leq n$, let
$\bm g_i\in\mathbb Q^{4n}$ equal $\bm g$ in its $i$th four-dimensional
block and zero elsewhere.  Hence
\begin{equation}
 \bm g_i^{\mathsf T}\bm g_j
 =L^2\mathbf1_{\{i=j\}}.
 \label{eq:svbc-orthogonal-tags}
\end{equation}
Use the orthogonal block decomposition
\begin{equation}
 \mathbb R^{d_0}
 =
 \underbrace{\mathbb R^n}_{\mathrm{identity}}
 \oplus
 \underbrace{\mathbb R^{3q}}_{\mathrm{incidence}}
 \oplus
 \underbrace{\mathbb R^{4n}}_{\mathrm{tag}}
 \oplus
 \underbrace{\mathbb R^4}_{\mathrm{height}},
 \qquad d_0=5n+3q+4,
 \label{eq:svbc-coordinate-blocks}
\end{equation}
and define, for $1\leq i\leq n$,
\begin{align}
 \bm p_i
 &=\bm e_i\oplus\bm b_i\oplus\bm0_{4n}\oplus\bm0_4,
 \label{eq:svbc-candidate}\\
 \bm w_i
 &=\bm0_n\oplus\bm b_i\oplus\bm g_i\oplus\bm0_4,
 \label{eq:svbc-paired-anchor}\\
 \bm a
 &=\bm0_n\oplus\bm u\oplus\bm0_{4n}\oplus\bm h.
 \label{eq:svbc-global-anchor}
\end{align}
Here $\bm e_i\in\mathbb R^n$ and $\bm b_i\in\{0,1\}^{3q}$ is the
$i$th incidence column.  
The face $\operatorname{conv}\{\bm p_1,\ldots,\bm p_n\}$ of the
simplex defined below will be called the candidate face.  Let
\begin{equation}
 \mathcal P
 =\operatorname{conv}
 \{\bm a,\bm p_1,\ldots,\bm p_n,\bm w_1,\ldots,\bm w_n\}.
 \label{eq:svbc-embedded-simplex}
\end{equation}

These $2n+1$ points are affinely independent, as seen by
projection onto the tag and identity blocks. Thus $\mathcal P$ is a $2n$-simplex in $\mathbb R^{d_0}$.
We shall prove that
\begin{equation}
 R^*(\mathcal P)\geq\rho_q
 \quad\Longleftrightarrow\quad
 (U,\mathcal S)\text{ has an exact cover}.
\label{eq:embedded-threshold-equivalence}
\end{equation}

\textbf{Completeness.}
We first show that an exact cover gives a threshold-feasible point
on the candidate face.  Let $C$ index an exact cover, and let
$\bm\lambda^C$ be the normalized indicator vector defined in
\cref{eq:l1-cover-point}.  Then
$\|\bm\lambda^C\|_2^2=1/q$ and $\bm B\bm\lambda^C=\bm u$.
We embed this vector into the candidate face by setting
\begin{equation}
 \bm x_C
 =\sum_{i=1}^n\lambda_i^C\bm p_i
 =\bm\lambda^C\oplus\bm u\oplus\bm0_{4n}\oplus\bm0_4.
 \label{eq:svbc-yes-point}
\end{equation}
This point belongs to $\mathcal P$.
Using these identities and the orthogonality of the coordinate blocks gives
\begin{align}
 &\|\bm x_C-\bm p_i\|_2^2
 =\rho_q+\frac2q\mathbf1_{\{i\notin C\}},
 \notag\\
 &\|\bm x_C-\bm w_i\|_2^2
 =\frac1q+3-\frac3q+L^2=\rho_q,\quad 
 \|\bm x_C-\bm a\|_2^2
 =\frac1q+H^2=\rho_q.
 \label{eq:svbc-completeness-distances}
\end{align}
Therefore $R^*(\mathcal P)\geq\rho_q$.

\textbf{Soundness.}
We prove the converse in two steps.  Weighted averages of squared
distances first show that the barycentric coefficients of $\bm a$
and $\bm w_i$ vanish.  The remaining inequalities then yield the
two norm equalities required by \cref{lem:exact-cover-rigidity}.
Let $\bm x\in\mathcal P$ satisfy
$\|\bm x-\bm v\|_2^2\geq\rho_q$ for every vertex $\bm v$ of
$\mathcal P$, and write its unique barycentric representation as
\begin{equation}
 \bm x=s\bm a+\sum_{i=1}^n\lambda_i\bm p_i
                 +\sum_{i=1}^n\mu_i\bm w_i,
 \qquad
 s\geq0,\quad\bm\lambda,\bm\mu\in\mathbb R_+^n,\quad
 s+\bm1_n^{\mathsf T}(\bm\lambda+\bm\mu)=1.
 \label{eq:svbc-barycentric}
\end{equation}
Define the coefficient sums $\Lambda,M,T$ and the auxiliary
quantities $E,V$ by
\begin{align}
 \Lambda&=\bm1_n^{\mathsf T}\bm\lambda,\qquad
 M=\bm1_n^{\mathsf T}\bm\mu,\qquad
 T=\Lambda+M=1-s,\notag\\
 E&=\|\bm\lambda\|_2^2+
    \|\bm B(\bm\lambda+\bm\mu)\|_2^2+
    L^2\|\bm\mu\|_2^2, \qquad
 V=3+\frac{\Lambda+ML^2}{T}\quad(T>0).
 \label{eq:svbc-energy-notation}
\end{align}
If $T=0$, then $\bm x=\bm a$ and
$\|\bm x-\bm a\|_2^2=0<\rho_q$, contradicting threshold feasibility.
Hence $T>0$.  Dividing the coefficients $\lambda_i,\mu_i$ by
their sum $T$ gives the convex combination
\begin{equation}
 \bm y=\frac1T\sum_{i=1}^n
          (\lambda_i\bm p_i+\mu_i\bm w_i)
       ,
 \qquad
 \bm x=s\bm a+T\bm y.
\label{eq:y}
\end{equation}
The column sums of $\bm B$ and the orthogonal tag blocks give
$\bm1_{3q}^{\mathsf T}\bm B(\bm\lambda+\bm\mu)=3T$ and
$\|\sum_i\mu_i\bm g_i\|_2^2=L^2\|\bm\mu\|_2^2$.
Together with $\|\bm u\|_2^2=3/q$, these identities yield
\[
 \|\bm x\|_2^2
 =E+\frac{6sT+3s^2}{q}+s^2H^2,
 \qquad
 \bm x^{\mathsf T}\bm a=\frac3q+sH^2.
\]
Expanding the distance to $\bm a$ and the weighted average of
distances to $\bm p_i,\bm w_i$, and using
$T\bm y=\bm x-s\bm a$, we obtain
\begin{align}
 D_a&:=\|\bm x-\bm a\|_2^2=E+T^2K,\notag\\
 \overline D
 &:=\frac1T\sum_{i=1}^n
   \bigl(\lambda_i\|\bm x-\bm p_i\|_2^2+
         \mu_i\|\bm x-\bm w_i\|_2^2\bigr)\notag
=-\left(\frac2T-1\right)E+(1-T)^2K+V.
 \label{eq:svbc-face-averages}
\end{align}
Since $0<T\leq1$, $(1-T/2)D_a+(T/2)\overline D$
is a convex combination of squared distances to the vertices.
The coefficients cancel $E$, giving
\begin{align}
 \min_{\bm v\in\operatorname{vert}(\mathcal P)}
       \|\bm x-\bm v\|_2^2
 \leq
 \left(1-\frac T2\right)D_a+\frac T2\overline D
  =\rho_q-s\rho_q-\frac Mq.
 \label{eq:svbc-face-forcing}
\end{align}
For a threshold-feasible point, the left-hand side of
\cref{eq:svbc-face-forcing} is at least $\rho_q$.
Since $s,M\geq0$ and $\rho_q,q>0$, this forces $s=M=0$.
Consequently, \(\Lambda=1\) and
\begin{equation}
 \bm x=\sum_{i=1}^n\lambda_i\bm p_i,\qquad
 \bm\lambda\in\Delta_{n-1}^{+}.
 \label{eq:svbc-candidate-face}
\end{equation}

We have now restricted every threshold-feasible point to the
candidate face.  It remains to recover the two rigidity equalities
from its barycentric coefficients.
Set $A=\|\bm\lambda\|_2^2$ and $C=\|\bm B\bm\lambda\|_2^2$.
The global-anchor constraint
$\|\bm x-\bm a\|_2^2\geq\rho_q$
and the $\lambda_i$-weighted average of the candidate constraints
$\|\bm x-\bm p_i\|_2^2\geq\rho_q$ give, respectively,
\begin{equation}
 \rho_q\leq A+C+K,\qquad
 \rho_q\leq
 \sum_{i=1}^n\lambda_i\|\bm x-\bm p_i\|_2^2
 =4-A-C.
 \label{eq:svbc-energy-sum}
\end{equation}
The two inequalities of \cref{eq:svbc-energy-sum} imply
\begin{equation}
 A+C=\frac4q.
 \label{eq:svbc-total-energy}
\end{equation}
The constraints $\|\bm x-\bm w_i\|_2^2\geq\rho_q$
for the paired anchors $\bm w_i$, give
\begin{equation}
 \rho_q
 \leq\sum_{i=1}^n\lambda_i\|\bm x-\bm w_i\|_2^2
 =A+3-C+L^2.
 \label{eq:svbc-paired-average}
\end{equation}
The Cauchy--Schwarz bound used in the proof of
\cref{lem:exact-cover-rigidity} gives $C\geq3/q$.
Substituting $A=4/q-C$ and the definitions of $L^2$ and $\rho_q$
into \cref{eq:svbc-paired-average} gives $C\leq3/q$.
Together with the lower bound, this yields
\begin{equation}
 \|\bm\lambda\|_2^2=\frac1q,\qquad
 \|\bm B\bm\lambda\|_2^2=\frac3q.
 \label{eq:svbc-separated-energies}
\end{equation}
By \cref{lem:exact-cover-rigidity}, the support of
$q\bm\lambda$ indexes an exact cover, proving
\cref{eq:embedded-threshold-equivalence}.
We next make the simplex full dimensional while preserving
this equivalence.

\subsection{Full-dimensional lifting via a balanced spherical cluster}\label{sec:full-dimensional}

We replace $\bm a$ by anchors
$\bm a_j=\bar{\bm a}+\bm c_j$ $(0\leq j\leq r)$
while preserving threshold feasibility.

\textbf{Lifting construction.}
Put
\begin{equation}
 \varepsilon=\frac1q,\qquad
 H_0^2=H^2-\varepsilon^2,\qquad
 \bar{\bm a}
 =\bm0_n\oplus\bm u\oplus\bm0_{4n}\oplus\bm h_0,
 \label{eq:svbc-lowered-anchor}
\end{equation}
where $\bm h_0\in\mathbb Q^4$ satisfies
$\|\bm h_0\|_2^2=H_0^2$. The rational construction of $\bm h_0$ is included in
\cref{lem:svbc-realization}. The decrease by \(\varepsilon^2\) reserves exactly the squared-distance
contribution of the orthogonal cluster vectors introduced below. 

Define the direction space of the preceding construction by
\begin{equation}
 \mathcal L
 =\operatorname{span}
 \{\bm p_i-\bar{\bm a},\bm w_i-\bar{\bm a}:1\leq i\leq n\},
 \qquad
 W_0=\mathcal L^\perp\subseteq\mathbb R^{d_0}.
 \label{eq:svbc-base-space}
\end{equation}

The same tag- and identity-block argument used in
\cref{sec:embedded-reduction} shows that the \(2n\) vectors
\[
 \{\bm p_i-\bar{\bm a},\bm w_i-\bar{\bm a}:1\leq i\leq n\}
\]
are linearly independent.  Hence
\begin{equation}
 \dim\mathcal L=2n,\qquad
 s_0:=\dim W_0=d_0-2n=3n+3q+4.
 \label{eq:svbc-base-dimensions}
\end{equation}

The rational cluster construction in \cref{app:svbc-realization}
uses nine additional coordinates to obtain a common norm and a
strictly positive balancing relation.  We therefore work in
$\mathbb R^d=\mathbb R^{d_0}\oplus\mathbb R^9$, where
$d=d_0+9=5n+3q+13$, and extend all base points by nine zero
coordinates.  We retain their symbols and likewise identify
$\mathcal L$ and $W_0$ with their zero extensions.
The orthogonal complement of $\mathcal L$ in the enlarged ambient space is
therefore
\begin{equation}
 W=\mathcal L^\perp=W_0\oplus\mathbb R^9,\qquad
 r:=\dim W=s_0+9=d-2n.
 \label{eq:svbc-lifted-dimension}
\end{equation}

By \cref{lem:svbc-realization}, proved in
\cref{app:svbc-realization}, one can construct in polynomial
time $r+1$ affinely independent vectors
$\bm c_0,\ldots,\bm c_r\in W\cap\mathbb Q^d$ such that
\begin{equation}
 \|\bm c_j\|_2^2=\varepsilon^2
 \quad(0\leq j\leq r),\qquad
 \bm0_d\in\operatorname{relint}
 \left(
   \operatorname{conv}\{\bm c_0,\ldots,\bm c_r\}
 \right).
 \label{eq:svbc-cluster-properties}
\end{equation}
We call the cluster $\{\bm c_0,\ldots,\bm c_r\}$ balanced
because there exist coefficients $\theta_j>0$ with
$\sum_{j=0}^r\theta_j=1$ and
$\sum_{j=0}^r\theta_j\bm c_j=\bm0_d$.

We remove the original global anchor $\bm a$ and replace it by $\bm a_j=\bar{\bm a}+\bm c_j , 0\leq j\leq r$. The candidate vertices $\bm p_i$ and paired anchors $\bm w_i$ are left
unchanged.  The resulting polytope is
\begin{equation}
 \mathcal Q
 =\operatorname{conv}
 \{\bm p_1,\ldots,\bm p_n,\bm w_1,\ldots,\bm w_n,
   \bm a_0,\ldots,\bm a_r\}.
 \label{eq:svbc-full-simplex}
\end{equation}

\textbf{Full dimensionality.}
We verify that the cluster supplies exactly the missing directions.
The balancing property gives
$\bar{\bm a}\in\operatorname{conv}\{\bm a_0,\ldots,\bm a_r\}
\subseteq\mathcal Q$.  Relative to this point, the differences
$\bm p_i-\bar{\bm a}$ and $\bm w_i-\bar{\bm a}$ span $\mathcal L$,
while $\bm a_j-\bar{\bm a}=\bm c_j$ span $W$.
Since $\mathcal L\oplus W=\mathbb R^d$, the affine hull of
$\mathcal Q$ is $\mathbb R^d$.  The construction has
$2n+(r+1)=d+1$ points, so they are affinely independent and
$\mathcal Q$ is a full-dimensional simplex.

We next prove that the lifting preserves the reduction:
\begin{equation}
 R^*(\mathcal Q)\geq\rho_q
 \quad\Longleftrightarrow\quad
 (U,\mathcal S)\text{ has an exact cover}.
 \label{eq:svbc-full-threshold-equivalence}
\end{equation}

\textbf{Completeness.}
Suppose that \(C\) indexes an exact cover.
We show that the point $\bm x_C$ from
\cref{eq:svbc-yes-point} remains feasible for $\mathcal Q$.
By construction,
\(\bm x_C\in\operatorname{conv}\{\bm p_1,\ldots,\bm p_n\}
\subseteq\mathcal Q\), and its squared distances to \(\bm p_i\)
and \(\bm w_i\) are still given by
\cref{eq:svbc-completeness-distances}. Moreover,
\[
 \bm x_C-\bar{\bm a}
 =\frac1q\sum_{i\in C}(\bm p_i-\bar{\bm a})
 \in\mathcal L.
\]
Since \(\bm c_j\in W=\mathcal L^\perp\), for every \(j\),
\begin{equation}
 \|\bm x_C-\bm a_j\|_2^2
 =\|\bm x_C-\bar{\bm a}\|_2^2+\|\bm c_j\|_2^2
 =\frac1q+H_0^2+\varepsilon^2
 =\rho_q.
 \label{eq:svbc-cluster-completeness}
\end{equation}
Thus \(R^*(\mathcal Q)\geq\rho_q\).

\textbf{Soundness.}
We show that a threshold-feasible point still lies in
$\operatorname{conv}\{\bm p_1,\ldots,\bm p_n\}$.
A weighted average of squared distances to the new anchors
$\bm a_j$ replaces the squared distance to $\bm a$ in the earlier
argument.
Suppose that $\bm x\in\mathcal Q$ is threshold feasible. Let
\begin{equation}
 \delta(\bm x)
 :=\min_{\bm v\in\operatorname{vert}(\mathcal Q)}
      \|\bm x-\bm v\|_2^2
 \geq\rho_q.
 \label{eq:svbc-lifted-threshold-feasibility}
\end{equation}
Since $\mathcal Q$ is a simplex, write $\bm x$ as
\begin{align}
 \bm x
 &=\sum_{i=1}^n\lambda_i\bm p_i+
   \sum_{i=1}^n\mu_i\bm w_i+
   \sum_{j=0}^r\alpha_j(\bar{\bm a}+\bm c_j),\notag\\
 &\lambda_i,\mu_i,\alpha_j\geq0,\qquad
 \sum_{i=1}^n(\lambda_i+\mu_i)+\sum_{j=0}^r\alpha_j=1.
 \label{eq:svbc-cluster-barycentric}
\end{align}
Set
\begin{equation}
 s=\sum_{j=0}^r\alpha_j,\qquad
 \Lambda=\sum_{i=1}^n\lambda_i,\qquad
 M=\sum_{i=1}^n\mu_i,\qquad
 T=\Lambda+M=1-s,
 \label{eq:svbc-lifted-masses}
\end{equation}
and let
\begin{equation}
 \bm z_\perp=\sum_{j=0}^r\alpha_j\bm c_j\in W.
 \label{eq:svbc-perpendicular-component}
\end{equation}

The triangle inequality and the common-norm property in
\cref{eq:svbc-cluster-properties} give
\begin{equation}
 \|\bm z_\perp\|_2
 \leq\sum_{j=0}^r\alpha_j\|\bm c_j\|_2
 =s\varepsilon,
 \qquad
 \|\bm z_\perp\|_2^2\leq s^2\varepsilon^2.
 \label{eq:svbc-perpendicular-bound}
\end{equation}

By the relative-interior property in
\cref{eq:svbc-cluster-properties}, there exist coefficients
$\theta_j>0$ satisfying
\begin{equation}
 \sum_{j=0}^r\theta_j=1,\qquad
 \sum_{j=0}^r\theta_j\bm c_j=\bm0_d.
 \label{eq:svbc-cluster-weights}
\end{equation}
The balancing identity $\sum_j\theta_j\bm c_j=\bm0_d$ cancels
the cross terms in the averaged squared distances.  Since
$\|\bm c_j\|_2^2=\varepsilon^2$ and $\sum_j\theta_j=1$, we obtain
\begin{equation}
 D_{\mathcal A}
 :=\sum_{j=0}^r\theta_j\|\bm x-\bm a_j\|_2^2
 =\|\bm x-\bar{\bm a}\|_2^2+\varepsilon^2.
 \label{eq:svbc-cluster-average}
\end{equation}
If $T=0$, then nonnegativity gives
$\bm\lambda=\bm\mu=\bm0_n$ and $s=1$, so
$\bm x-\bar{\bm a}=\bm z_\perp$.  By
\cref{eq:svbc-perpendicular-bound,eq:svbc-cluster-average},
\begin{equation}
 \delta(\bm x)
 \leq D_{\mathcal A}
 \leq2\varepsilon^2
 =\frac2{q^2}
 <\rho_q,
\end{equation}
contradicting \cref{eq:svbc-lifted-threshold-feasibility}.  Hence $T>0$.
Retain $E$, $V$ from \cref{eq:svbc-energy-notation} and $\bm y$ from \cref{eq:y}.
Set
\[
 K_0:=H_0^2-\frac3q=K-\varepsilon^2.
\]
The barycentric representation in \cref{eq:y} becomes
\[
 \bm x=T\bm y+s\bar{\bm a}+\bm z_\perp,
 \qquad
 \bm x-\bar{\bm a}
 =
 T(\bm y-\bar{\bm a})+\bm z_\perp.
\]
Here $\bm y-\bar{\bm a}\in\mathcal L$ and
$\bm z_\perp\in W=\mathcal L^\perp$, so their squared norms add.
Likewise, $T\bm y+s\bar{\bm a}-\bm v\in\mathcal L$ for every
$\bm v\in\{\bm p_i,\bm w_i:1\leq i\leq n\}$, giving
$\|\bm x-\bm v\|_2^2
 =\|T\bm y+s\bar{\bm a}-\bm v\|_2^2+\|\bm z_\perp\|_2^2$.
Applying \cref{eq:svbc-cluster-average,eq:svbc-face-averages}
with $\bm a$ replaced by $\bar{\bm a}$ and $K$ by $K_0$ therefore yields
\begin{align*}
 D_{\mathcal A}
 &:=\sum_{j=0}^r\theta_j\|\bm x-\bm a_j\|_2^2
 =\|\bm x-\bar{\bm a}\|_2^2+\varepsilon^2=E+T^2K_0+\|\bm z_\perp\|_2^2+\varepsilon^2,\\
 D_{\mathcal N}
 &:=\frac1T\sum_{i=1}^n
   \bigl(\lambda_i\|\bm x-\bm p_i\|_2^2+
         \mu_i\|\bm x-\bm w_i\|_2^2\bigr)=-\left(\frac2T-1\right)E+s^2K_0+V+\|\bm z_\perp\|_2^2.
 \label{eq:svbc-lifted-averages}
\end{align*}

Both \(D_{\mathcal A}\) and \(D_{\mathcal N}\) are convex combinations
of vertex distances.  Since \(0<T\leq1\), the coefficients
\(1-T/2\) and \(T/2\) also form a convex combination; they are chosen
to cancel \(E\).  Using the same cancellation as in
\cref{eq:svbc-face-forcing} and
\(K_0=K-\varepsilon^2\), we obtain
\begin{align}
 \delta(\bm x)&\leq
 \left(1-\frac T2\right)D_{\mathcal A}
 +\frac T2D_{\mathcal N}\notag
 =\|\bm z_\perp\|_2^2
   +T\rho_q-\frac Mq+s\varepsilon^2\notag\\
 &\leq
 \rho_q-\frac Mq
 -s\bigl[\rho_q-(1+s)\varepsilon^2\bigr],
 \label{eq:svbc-lifted-face-forcing}
\end{align}
where the last inequality follows from
\(\|\bm z_\perp\|_2^2\leq s^2\varepsilon^2\) and \(T=1-s\).

Since $0\leq s\leq1$ and $q\geq3$,
\begin{equation}
 \rho_q-(1+s)\varepsilon^2
 \geq\rho_q-\frac2{q^2}>0.
 \label{eq:svbc-lifted-positive-gap}
\end{equation}
Threshold feasibility gives $\delta(\bm x)\geq\rho_q$.
Thus \cref{eq:svbc-lifted-face-forcing,eq:svbc-lifted-positive-gap}
force $s=M=0$ and $\Lambda=1$.  Nonnegativity then gives
$\alpha_j=\mu_i=0$ for all $i,j$, so
$\bm x=\sum_i\lambda_i\bm p_i$ with
$\bm\lambda\in\Delta_{n-1}^{+}$ lies on the original candidate face.

For such a point, $\bm x-\bar{\bm a}\in\mathcal L$.
By orthogonality and $H_0^2+\varepsilon^2=H^2$, we have
\begin{equation}
 \|\bm x-\bm a_j\|_2^2
 =
 \|\bm x-\bar{\bm a}\|_2^2+\varepsilon^2
 =
 \|\bm x-\bm a\|_2^2
 \qquad(0\leq j\leq r),
 \label{eq:svbc-restored-anchor}
\end{equation}
where the original anchor $\bm a$ is identified with its zero
extension.  The distances to $\bm p_i$ and $\bm w_i$ are unchanged.
Thus, after removing the added zero coordinates, $\bm x$ is a
threshold-feasible point of the embedded simplex $\mathcal P$.
By \cref{eq:embedded-threshold-equivalence}, the X3C instance
has an exact cover.  This proves
\cref{eq:svbc-full-threshold-equivalence}.  The geometric reduction
is now complete; it remains to justify the rational construction
and the numerical bounds needed for strong NP-hardness.

\subsection{Rational realization and integer encoding}
\label{sec:svbc-encoding}

To obtain integer data, we need a rational realization with
a polynomially bounded common denominator.
The following lemma, proved in \cref{app:svbc-realization},
provides one.

\begin{lemma}[Rational realization]
\label{lem:svbc-realization}
For the X3C instance above, the vectors
$\bm g,\bm h,\bm h_0\in\mathbb Q^4$ can be constructed in polynomial
time with the squared norms prescribed in
\cref{eq:svbc-parameters,eq:svbc-lowered-anchor},
with $q\bm g,q\bm h,q\bm h_0$ integral and with entries of
$\bm g,\bm h,\bm h_0$ bounded by an absolute constant.
For the resulting space $W$, one can construct in polynomial time
affinely independent vectors $\bm c_0,\ldots,\bm c_r\in W\cap\mathbb Q^d$
satisfying \cref{eq:svbc-cluster-properties}, together with a
positive integer $\Gamma$ such that
\[
 q\Gamma\bm c_j\in\mathbb Z^d\quad(0\leq j\leq r),
 \qquad
 \Gamma=O((n+q)^2q^3).
\]
\end{lemma}

We now apply these bounds to the final scaling.
Scale every vertex of $\mathcal Q$ by $q\Gamma$ and the
threshold $\rho_q$ by $(q\Gamma)^2$.
Since $q\bm p_i$, $q\bm w_i$, $q\bar{\bm a}$, and
$q\Gamma\bm c_j$ are integral, this produces a full-dimensional
integer simplex.
Squared distances scale by $(q\Gamma)^2$, preserving
\cref{eq:svbc-full-threshold-equivalence}.

The vectors $\bm p_i,\bm w_i,\bar{\bm a}$ have bounded
coordinates, and $\|\bm c_j\|_2=1/q$, so
\[
 \max_{\bm v\in\operatorname{vert}(\mathcal Q)}
 \|q\Gamma\bm v\|_\infty=O(q\Gamma),
 \qquad
 (q\Gamma)^2\rho_q=4q(q-1)\Gamma^2.
\]
Since $\Gamma=O((n+q)^2q^3)$ and $d=5n+3q+13$, the number
$d(d+1)$ of output coordinates, their magnitudes, and the integer
threshold $(q\Gamma)^2\rho_q$ are polynomially bounded in the
explicit X3C input size.
This proves strong NP-hardness for full-dimensional integer
simplices.

We next establish exact algorithms and short certificates,
completing the two NP-completeness proofs.

\section{Exact solution and NP membership over the unit
\texorpdfstring{$\ell_1$}{l1}-ball}
\label{sec:l1-exact}

In this section, we give an exact algorithm for $\ell_1$-WMD and prove that
$\ell_1$-WMD-DEC belongs to $\mathrm{NP}$.

Let $\ell_{\mathrm{WMD}}$ denote the total binary
encoding length, including the dimensions, sites, positive weights, and,
for the decision problem, the threshold.

\begin{theorem}
\label{thm:l1-exact-algorithm}
The problem $\ell_1$-WMD can be solved exactly in
$2^d(r+d)^{O(d)}\cdot\operatorname{poly}(\ell_{\mathrm{WMD}})$
bit operations. There exists an optimal solution of encoding length
polynomial in $\ell_{\mathrm{WMD}}$ whose coordinates lie in a
common real algebraic extension of $\mathbb Q$ of degree at most two.
\end{theorem}

In particular, $\ell_1$-WMD is solvable in polynomial time
for every fixed $d$.

\subsection{Orthant decomposition and exact quadratic lifting}
\label{sec:l1-exact-lifting}

For each sign vector $\bm\sigma\in\{-1,1\}^d$, define the orthant
simplex
\begin{equation}
 P_{\bm\sigma}
 :=
 \left\{
   \bm x\in\mathbb R^d:
   \sigma_i x_i\geq0\ (1\leq i\leq d),\quad
   \bm\sigma^{\mathsf T}\bm x\leq1
 \right\}.
 \label{eq:l1-exact-orthant-simplex}
\end{equation}
If $\bm x\in P_{\bm\sigma}$, then $\|\bm x\|_1
 =
 \sum_{i=1}^d\sigma_i x_i
 =
 \bm\sigma^{\mathsf T}\bm x
 \leq1 $.
Conversely, every $\bm x\in B_1^d$ belongs to at least one
$P_{\bm\sigma}$, obtained by choosing $\sigma_i x_i=|x_i|$.
Consequently,
\begin{equation}
 B_1^d
 =
 \bigcup_{\bm\sigma\in\{-1,1\}^d}P_{\bm\sigma}.
 \label{eq:l1-exact-orthant-union}
\end{equation}
To separate the quadratic term from the distance inequalities,
introduce $z=\|\bm x\|_2^2$ and a variable $t$ to be maximized subject
to $t\leq\omega_s\|\bm x-\bm a_s\|_2^2$ for every site.
For fixed $\bm x$, the largest admissible $t$ is the original
dispersion objective.
For every site, define the affine function $ L_s(\bm x,z)
 :=
 \omega_s z
 -2\omega_s\bm a_s^{\mathsf T}\bm x
 +\omega_s\|\bm a_s\|_2^2$.
On the paraboloid $z=\|\bm x\|_2^2$, this gives
$L_s(\bm x,z)=\omega_s\|\bm x-\bm a_s\|_2^2$,
so the lifting preserves the distance constraints exactly.

For fixed $\bm\sigma$, define the rational polytope
\begin{equation}
 \begin{aligned}
 K_{\bm\sigma}:=\bigl\{(\bm x,z,t)\in\mathbb R^{d+2}:\
   &\bm x\in P_{\bm\sigma},\quad 0\leq z\leq1,\\
   &0\leq t\leq L_s(\bm x,z)\quad(1\leq s\leq r)\bigr\}.
 \end{aligned}
 \label{eq:l1-exact-lifted-polyhedron}
\end{equation}
The equation $z=\|\bm x\|_2^2$ is represented separately by
\begin{equation}
 Q(\bm x,z,t)
 :=
 z-\|\bm x\|_2^2
 =
 0.
 \label{eq:l1-exact-paraboloid}
\end{equation}
The bounds $0\leq z\leq1$ hold automatically on the paraboloid
above $P_{\bm\sigma}$, and $t\geq0$ preserves the optimal value
because all weighted squared distances are nonnegative.
Moreover, $K_{\bm\sigma}$ is bounded: $\bm x$ and $z$ are bounded,
and $0\leq t\leq L_1(\bm x,z)$ bounds $t$ as well.
Using \cref{eq:l1-exact-orthant-union} and the distance identity above,
\begin{equation}
 \Phi^*(\mathcal A,\bm\omega)
 =
 \max_{\bm\sigma\in\{-1,1\}^d}
 \max
 \left\{
   t:
   (\bm x,z,t)\in K_{\bm\sigma},\
   Q(\bm x,z,t)=0
 \right\}.
 \label{eq:l1-exact-lifted-equivalence}
\end{equation}

Thus $\ell_1$-WMD reduces to $2^d$ subproblems, one for each
$\bm\sigma\in\{-1,1\}^d$, in which we maximize $t$ over
the polytope $K_{\bm\sigma}$ subject to $Q(\bm x,z,t)=0$.

\subsection{Linear optimization over a polytope and one quadric}
\label{sec:one-quadric-optimization}

For the orthant subproblem indexed by $\bm\sigma$, we apply
the following theorem with $P=K_{\bm\sigma}$,
$\bm y=(\bm x,z,t)$, $q=Q$, and $\ell(\bm y)=t$.

\begin{theorem}[One-quadric optimization]
\label{thm:one-quadric-optimization}
Let $P=\{\bm y\in\mathbb R^N:\bm C\bm y\leq\bm b\} $
be a nonempty rational polytope described by $M$ inequalities, and let
\begin{equation}
 q(\bm y)=\gamma+2\bm g^{\mathsf T}\bm y
          -\bm y^{\mathsf T}\bm H\bm y,
 \qquad \bm H\succeq\bm0_N,
 \label{eq:one-quadric-general-form}
\end{equation}
where $\bm H,\bm g,\gamma$ are rational.  Suppose that
$S:=P\cap\{q=0\}$ is nonempty.  For every rational affine function
$\ell(\bm y)=\ell_0+\bm c^{\mathsf T}\bm y$, the problem
\begin{equation}
 \max\{\ell(\bm y):\bm y\in S\}
 \label{eq:one-quadric-linear-optimization}
\end{equation}
can be solved exactly in $M^{O(N)}\operatorname{poly}(L)$ bit
operations, where $L$ is the total binary encoding length of the data.
Moreover, an optimal solution can be represented in one common real
algebraic extension of degree at most two, with encoding length
polynomial in $L$.
\end{theorem}

\begin{proof}
We first locate an optimizer in the relative interior of a face of
$P$, so that the active linear constraints can be absorbed into an
affine parametrization.  We then solve the resulting regular,
singular, and constant-objective branches.  Their feasible candidates
include an optimizer, and exact comparison selects a largest value.

\textbf{Step 1: locating the minimal face.}
Since $S$ is nonempty and compact, choose a global maximizer
$\bm y^*$. Let $F^*$ be the minimal face of $P$ containing
$\bm y^*$. We show that its affine hull $E^*=\operatorname{aff}F^*$
is determined by the constraints active at $\bm y^*$.

Define the active index set
\[
 \mathcal A^*:=\{i\in\{1,\ldots,M\}:\bm C_i\bm y^*=b_i\}.
\]
The active equalities give
\[
 F^*=\{\bm y\in P:\bm C_i\bm y=b_i\ (i\in\mathcal A^*)\},
 \qquad
 E^*=\{\bm y\in\mathbb R^N:\bm C_i\bm y=b_i\
                           (i\in\mathcal A^*)\}.
\]
All inequalities with $i\notin\mathcal A^*$ are strict at
$\bm y^*$. Hence every point of $E^*$ sufficiently close to
$\bm y^*$ lies in $P$. Consequently, $\bm y^*$ is a local
maximizer of $\ell$ on $E^*\cap\{q=0\}$.

Choose $J\subseteq\mathcal A^*$ whose constraint rows form a basis
of the row space of $\bm C_{\mathcal A^*}$.  Consistency at
$\bm y^*$ gives
\begin{equation}
 |J|\leq N,\qquad
 E^*=\{\bm y:\bm C_J\bm y=\bm b_J\}.
 \label{eq:one-quadric-face-basis}
\end{equation}
We therefore enumerate all linearly independent sets of at most $N$
constraint rows, including the empty set.  Their number is at most $ \sum_{k=0}^{\min\{N,M\}}\binom{M}{k}=M^{O(N)}$.
The enumeration thus includes a set $J$ defining $E^*$.
We now describe the calculation for each enumerated set.

\textbf{Step 2: rational parametrization and kernel decomposition.}
Fix $J$ and write $E_J=\{\bm y:\bm C_J\bm y=\bm b_J\}$.
If this system is inconsistent, discard $J$.  Otherwise, rational
Gaussian elimination gives a particular solution $\bm y_0$ and a
matrix $\bm V$ whose columns form a basis of $\ker\bm C_J$.  Thus
\begin{equation}
 \bm y=\bm y_0+\bm V\bm u,\qquad
 \bm u\in\mathbb R^p,\qquad p=N-\operatorname{rank}\bm C_J.
 \label{eq:one-quadric-affine-parametrization}
\end{equation}
Every point in $E_J$ can be written uniquely as
$\bm y=\bm y_0+\bm V\bm u$ with $\bm u\in\mathbb R^p$.  Substituting into the linear constraints
$\bm C\bm y\leq\bm b$, the quadratic equality $q(\bm y)=0$, and
the objective $\ell(\bm y)=\ell_0+\bm c^{\mathsf T}\bm y$ gives
\begin{align}
 P_J&:=\{\bm u:\bm C\bm V\bm u\leq\bm b-\bm C\bm y_0\},
 \label{eq:one-quadric-restricted-polyhedron}\\
 q_J(\bm u)&:=q(\bm y_0+\bm V\bm u)
       =\gamma_J+2\bm g_J^{\mathsf T}\bm u
                    -\bm u^{\mathsf T}\bm H_J\bm u,
 \label{eq:one-quadric-restricted-q}\\
 \ell_J(\bm u)&:=\ell(\bm y_0+\bm V\bm u)
       =\beta_J+\bm c_J^{\mathsf T}\bm u,
 \label{eq:one-quadric-restricted-objective}
\end{align}
where
\[
 \begin{aligned}
 \gamma_J&=\gamma+2\bm g^{\mathsf T}\bm y_0
                       -\bm y_0^{\mathsf T}\bm H\bm y_0,
 &\bm g_J&=\bm V^{\mathsf T}(\bm g-\bm H\bm y_0),\\
 \bm H_J&=\bm V^{\mathsf T}\bm H\bm V,
 &\bm c_J&=\bm V^{\mathsf T}\bm c,\qquad
 \beta_J=\ell_0+\bm c^{\mathsf T}\bm y_0.
 \end{aligned}
\]
For this $J$, the restricted optimization problem is
\[
 \max\bigl\{\ell_J(\bm u):
       \bm u\in P_J,\ q_J(\bm u)=0\bigr\}.
\]

Check whether the linear system
$\bm C\bm V\bm u\leq\bm b-\bm C\bm y_0$ is feasible.
If it is infeasible, proceed to the next $J$.
If $p=0$, retain $\bm y_0$ precisely when $q(\bm y_0)=0$,
then proceed to the next $J$.

In the remaining cases, $P_J\neq\varnothing$ and $p\geq1$.
The affine bijection $\bm u\mapsto\bm y_0+\bm V\bm u$ maps
$P_J$ onto $P\cap E_J$, so $P_J$ is a polytope.

We choose coordinates so that the quadratic term in $q_J$ depends
only on a block of variables $\bm u_R$.
Since $\bm H_J$ is rational and positive semidefinite, an
invertible rational change of coordinates gives the following
form. For simplicity, we retain the same notation after this change
of coordinates, including $\bm V$ in the parametrization.
\begin{equation}
 \bm H_J=
 \begin{bmatrix}
  \bm H_R&\bm0\\
  \bm0&\bm0
 \end{bmatrix},
 \qquad \bm H_R\succ\bm0.
 \label{eq:one-quadric-kernel-decomposition}
\end{equation}
Write $\bm u=(\bm u_R,\bm u_N)$, where
$\bm u_R\in\mathbb R^r$, $r=\operatorname{rank}\bm H_J$, and
$\bm u_N\in\mathbb R^{p-r}$ gives the coordinates along
$\ker\bm H_J$.
With the corresponding partitions
$\bm g_J=(\bm g_R,\bm g_N)$ and
$\bm c_J=(\bm c_R,\bm c_N)$, we obtain
\begin{align}
 q_J(\bm u_R,\bm u_N)
 &=\gamma_J+2\bm g_R^{\mathsf T}\bm u_R
             +2\bm g_N^{\mathsf T}\bm u_N
             -\bm u_R^{\mathsf T}\bm H_R\bm u_R,
 \label{eq:one-quadric-block-form}\\
 \ell_J(\bm u_R,\bm u_N)
 &=\beta_J+\bm c_R^{\mathsf T}\bm u_R
          +\bm c_N^{\mathsf T}\bm u_N.
 \label{eq:one-quadric-objective-block-form}
\end{align}
If $r=0$, the $R$-block and its terms are absent.
We now use these expressions to examine the regular, singular,
and constant-objective branches.

\textbf{Step 3: generating feasible candidates.}
For a set $J$ defining $E^*$, Step~1 shows that the corresponding
optimizer is a local maximizer of $\ell_J$ on $q_J=0$.
When $\bm c_J\neq\bm0$, we therefore first find stationary
points of $\ell_J$ on $q_J=0$, including singular points,
and then seek a point in $P_J$ from each resulting branch.
When $\bm c_J=\bm0$, we only need to determine whether
$P_J\cap\{q_J=0\}$ is nonempty.

\emph{Regular branches.}
Suppose that $\bm c_J\neq\bm0$. To find the regular stationary points of $\ell_J$ on $q_J=0$,
we solve the following Lagrange multiplier system:
\begin{equation}
 \bm c_J=2\lambda(\bm g_J-\bm H_J\bm u),
 \qquad q_J(\bm u)=0.
 \label{eq:one-quadric-stationarity}
\end{equation}
Here $\lambda\neq0$, and the block equations are
\begin{equation}
 \bm c_R=2\lambda(\bm g_R-\bm H_R\bm u_R),
 \qquad \bm c_N=2\lambda\bm g_N.
 \label{eq:one-quadric-block-stationarity}
\end{equation}

If $\bm g_N\neq\bm0$, the only possible multiplier is $\lambda=\frac{\bm g_N^{\mathsf T}\bm c_N}{2\|\bm g_N\|_2^2}$.
Discard the branch if $\lambda=0$ or
$\bm c_N\neq2\lambda\bm g_N$.  Otherwise,
\[
 \bm u_R=\bm H_R^{-1}
             \left(\bm g_R-\frac{\bm c_R}{2\lambda}\right)
\]
is fixed and rational.  Substitution into $q_J=0$ leaves the
linear equation
\begin{equation}
 2\bm g_N^{\mathsf T}\bm u_N
   =\bm u_R^{\mathsf T}\bm H_R\bm u_R
      -2\bm g_R^{\mathsf T}\bm u_R-\gamma_J.
 \label{eq:one-quadric-null-linear-equation}
\end{equation}
Solve this equation jointly with the inequalities of $P_J$ for
$\bm u_N$, retaining one feasible solution if any exists.
Since $\bm c_N=2\lambda\bm g_N$, all such solutions have the
same objective value.

If $\bm g_N=\bm0$, the second equation of \cref{eq:one-quadric-block-stationarity} requires
$\bm c_N=\bm0$; otherwise discard the branch.  Define
\begin{equation}
 C:=\gamma_J+\bm g_R^{\mathsf T}\bm H_R^{-1}\bm g_R,
 \qquad D:=\bm c_R^{\mathsf T}\bm H_R^{-1}\bm c_R>0.
 \label{eq:one-quadric-CD}
\end{equation}
Solving the first equation of \cref{eq:one-quadric-block-stationarity} and substituting the result into $q_J=0$ gives,
\begin{equation}
 \bm u_R=\bm H_R^{-1}\bm g_R
             -\frac1{2\lambda}\bm H_R^{-1}\bm c_R,
 \qquad C-\frac{D}{4\lambda^2}=0.
 \label{eq:one-quadric-regular-quadratic}
\end{equation}
There is no regular candidate when $C\leq0$.  When $C>0$,
the two possible multipliers are
\begin{equation}
 \lambda=\pm\frac12\sqrt{\frac DC}.
 \label{eq:one-quadric-lambda-roots}
\end{equation}
For each multiplier, $q_J$ and $\ell_J$ are independent of
$\bm u_N$, so we solve the linear inequalities
$(\bm u_R,\bm u_N)\in P_J$ for $\bm u_N$ and retain one feasible
point, if any.
A feasible point in the same quadratic extension can be found
by vertex enumeration, using
\cref{app:quadratic-arithmetic} to check the inequalities exactly.

\emph{Singular branch.}
Singular points satisfy $\nabla q_J=\bm0$.
If $\bm g_N\neq\bm0$, discard this branch.
Otherwise, $\bm u_R=\bm u_R^0:=\bm H_R^{-1}\bm g_R$, and
\[
 q_J(\bm u_R^0,\bm u_N)
 =C:=\gamma_J+\bm g_R^{\mathsf T}\bm H_R^{-1}\bm g_R.
\]
If $C\neq0$, discard the branch.
If $C=0$, solve the rational linear program
\[
 \max_{\bm u_N}
 \left\{\ell_J(\bm u_R^0,\bm u_N):
        (\bm u_R^0,\bm u_N)\in P_J\right\}.
\]
If feasible, retain an optimal vertex; otherwise discard the branch.

\emph{Constant-objective branch.}
Finally, if $\bm c_J=\bm0$, every feasible point has value $\beta_J$.
We therefore need only decide whether $P_J\cap\{q_J=0\}$ is
nonempty and, if so, construct one point in this intersection.
Set
\[
 h_{\min}:=\min_{\bm u\in P_J}q_J(\bm u),
 \qquad h_{\max}:=\max_{\bm u\in P_J}q_J(\bm u).
\]
Both extrema are attained because $P_J$ is compact.  Continuity
along segments in $P_J$ gives
\begin{equation}
 P_J\cap\{q_J=0\}\neq\varnothing
 \quad\Longleftrightarrow\quad h_{\min}\leq0\leq h_{\max}.
 \label{eq:one-quadric-interval-test}
\end{equation}
\Cref{app:quadric-feasibility} computes rational minimizers and
maximizers of $q_J$ on $P_J$ by vertex and active-set enumeration.
If either extremum is zero, its extremizer is a candidate;
if the extrema have opposite signs, solving $q_J=0$ on the
segment between the extremizers gives a candidate. This feasibility calculation
takes $M^{O(p)}\operatorname{poly}(L)$ bit operations.

For the enumerated set $J$ with $E_J=E^*$, the optimizer
$\bm u^*$ lies in $\operatorname{int}P_J$ and is covered by
the regular, singular, or constant-objective branch, which
retains a feasible point with the same objective value.
Thus comparing the retained candidates over all enumerated
sets $J$ yields a global optimum.

\textbf{Step 4: exact comparison and complexity.}
We now compare the objective values of the retained candidates.
Only the regular subcase $\bm g_N=\bm0$ can produce irrational
objective values.  A constant-objective candidate may have
irrational coordinates, but its value is $\beta_J\in\mathbb Q$;
the other branches admit rational candidates.

Each candidate value has the form $a+b\sqrt{\Delta}$ with
$a,b,\Delta\in\mathbb Q$ and $\Delta\geq0$.
Even when the radicands differ, these values can be compared
exactly in polynomial time using rational arithmetic, by checking
signs before squaring as described in
\cref{app:quadratic-arithmetic}.

Standard bit bounds for rational linear algebra and linear
programming, together with the quadratic formula, ensure that
the candidates can be chosen with polynomial encoding length.
The outer enumeration has size $M^{O(N)}$, and each branch uses
at most $M^{O(p)}\operatorname{poly}(L)$ bit operations, where
$p\leq N$.  Comparing all retained values and selecting a largest
one therefore takes $M^{O(N)}\operatorname{poly}(L)$ bit operations.
The selected candidate retains its representation in a common
real algebraic extension of $\mathbb Q$ of degree at most two.
\end{proof}

\subsection{Application to weighted
\texorpdfstring{$\ell_1$}{l1}-ball dispersion}
\label{sec:l1-exact-application}

\begin{proof}[Proof of \cref{thm:l1-exact-algorithm}]
For every $\bm\sigma\in\{-1,1\}^d$, $K_{\bm\sigma}$ is a rational
polytope and
\[
 (\bm0_d,0,0)\in K_{\bm\sigma}\cap\{Q=0\},
\]
since $L_s(\bm0_d,0)=\omega_s\|\bm a_s\|_2^2\geq0$.
With $\bm y=(\bm x,z,t)$, the quadratic equality is
\begin{equation}
 Q(\bm y)=2\bm g^{\mathsf T}\bm y-\bm y^{\mathsf T}\bm H\bm y,
 \qquad
 \bm H=\operatorname{diag}(\bm I_d,0,0)\succeq\bm0,
 \qquad \bm g=(\bm0_d,1/2,0).
 \label{eq:l1-exact-paraboloid-matrix}
\end{equation}
Thus \cref{thm:one-quadric-optimization} applies with
$P=K_{\bm\sigma}$, $q=Q$, and $\ell(\bm x,z,t)=t$.

The lifted dimension is $N=d+2$, and the number of linear
inequalities is $M=r+d+4$, including the two bounds on $z$.
Their coefficients have encoding length polynomial in
$\ell_{\mathrm{WMD}}$.  Solving all $2^d$ orthant problems and
comparing their values gives the global optimum by
\cref{eq:l1-exact-lifted-equivalence}, with total bit complexity
$2^d(r+d)^{O(d)}\operatorname{poly}(\ell_{\mathrm{WMD}})$.
The chosen optimal point belongs to one quadratic extension and
has polynomial encoding length by
\cref{thm:one-quadric-optimization}.  Its $\bm x$-coordinates give
the required solution of $\ell_1$-WMD.
\end{proof}

\subsection{NP membership}
\label{sec:l1-exact-np}

\begin{theorem}[NP membership]
\label{thm:l1-exact-np}
The decision problem $\ell_1$-WMD-DEC belongs to\/ $\mathrm{NP}$.
\end{theorem}

\begin{proof}
We describe a polynomial-length certificate for a point $\bm x$
satisfying $\|\bm x\|_1\leq1$ and
$\omega_s\|\bm x-\bm a_s\|_2^2\geq\tau$ for every $s=1,\ldots,r$,
and show that these inequalities can be verified exactly
in polynomial time.
Consider a YES-instance with threshold $\tau$.  By
\cref{thm:l1-exact-algorithm}, there exists an optimal point
$\bm x$ of polynomial encoding length whose coordinates belong to one
common real algebraic extension of degree at most two.  Hence every
coordinate can be represented as
\[
 x_i
 =
 \frac{p_i+q_i\sqrt{\Delta}}{r_i},
 \qquad
 p_i,q_i,r_i,\Delta\in\mathbb Z,
 \qquad
 r_i>0,\quad \Delta\geq0,
\]
where the same $\Delta$ is used for all coordinates.  A certificate is 
\(
 \left(
   \bm\sigma,\Delta,
   (p_i,q_i,r_i)_{i=1}^d
 \right),
\)
where $\bm\sigma\in\{-1,1\}^d$ specifies an orthant containing $\bm x$.

The verifier checks that $\Delta\geq0$, each $r_i>0$, and
$\bm\sigma\in\{-1,1\}^d$, and then checks
\[
 \sigma_i x_i\geq0
 \quad(1\leq i\leq d),
 \qquad
 \sum_{i=1}^d\sigma_i x_i\leq1.
\]
The sign checks ensure that
$\sum_i\sigma_i x_i=\sum_i|x_i|=\|\bm x\|_1$, so these conditions
verify membership in the unit $\ell_1$-ball.
For every site $\bm a_s$, the verifier also checks $ \omega_s\|\bm x-\bm a_s\|_2^2\geq\tau$.
Each check reduces to determining the sign of
$A+B\sqrt{\Delta}$ with rational $A$ and $B$, using the
single-radical sign test in \cref{app:quadratic-arithmetic}.
The number of arithmetic operations and the encoding lengths
of all intermediate rational quantities are polynomial in
the input and certificate lengths.  Hence the verifier runs
in polynomial time.

Every YES-instance has such a certificate by
\cref{thm:l1-exact-algorithm}.  Conversely, every accepted certificate
represents a point in the unit $\ell_1$-ball whose weighted squared
distance from every site is at least $\tau$, and hence certifies a
YES-instance.  Therefore $\ell_1$-WMD-DEC belongs to\/ $\mathrm{NP}$.
\end{proof}

Together with the strong NP-hardness established in
\cref{sec:l1-wmd}, this proves \cref{thm:l1-main}.

\section{Voronoi decomposition, exact solution, and NP membership}
\label{sec:voronoi}

Following the geometric argument of Casado et al.\
\cite[section~3]{CasadoEtAl2011}, we maximize the objective
over the vertices of the simplex intersected with each
Voronoi cell.
This yields an explicit candidate bound and polynomial-length
rational certificates for NP membership.

Let $\ell_{\mathrm{SVBC}}$ denote the total binary encoding length
of an instance, including the vertices, dimensions, and, for the
decision problem, the threshold.

\begin{theorem}
\label{thm:algorithmic}
For every rational $m$-simplex $S\subseteq\mathbb R^d$, the value
$R^*(S)$ is attained at a rational point of encoding length
polynomial in $\ell_{\mathrm{SVBC}}$.  Consequently, SVBC-DEC
belongs to $\mathrm{NP}$.  Moreover, SVBC can be solved exactly by
examining at most
\[
 N_m:=(m+1)\binom{2m+1}{m}
\]
rational candidate points, with total bit complexity $O\!\left(N_m\operatorname{poly}(d,\ell_{\mathrm{SVBC}})\right) $.
\end{theorem}

\begin{proof}
We first describe the Voronoi cells $\mathcal C_i$ in barycentric
coordinates $\bm\lambda$.  We then show that an optimum is attained
at a vertex of some $\mathcal C_i$ and enumerate these vertices
by solving linear systems.  Finally, we use a rational optimizer
$\bm\lambda$ of polynomial encoding length as a certificate
for NP membership.

\textbf{Step 1: barycentric Voronoi decomposition.}
Let
\[
 \bm V=[\bm v_0,\ldots,\bm v_m]\in\mathbb Q^{d\times(m+1)},
 \qquad
 \Delta_m=\{\bm\lambda\in\mathbb R_+^{m+1}:
                     \bm1_{m+1}^{\mathsf T}\bm\lambda=1\}.
\]
Since $S$ is a nondegenerate simplex, every point of $S$ has a
unique representation $\bm x=\bm V\bm\lambda$ with
$\bm\lambda\in\Delta_m$.  For each vertex index $i$, define its
nearest-vertex cell by
\begin{equation}
 \mathcal C_i
 =\{\bm\lambda\in\Delta_m:
       \|\bm V\bm\lambda-\bm v_i\|_2^2
       \leq\|\bm V\bm\lambda-\bm v_j\|_2^2
       \text{ for every }j\}.
 \label{eq:svbc-voronoi-cell}
\end{equation}
Expanding the distances cancels their common quadratic term:
\begin{equation}
 \mathcal C_i
 =\left\{\bm\lambda\in\Delta_m:
    2(\bm v_j-\bm v_i)^{\mathsf T}\bm V\bm\lambda
       \leq\|\bm v_j\|_2^2-\|\bm v_i\|_2^2
       \quad(j\neq i)\right\}.
 \label{eq:svbc-voronoi-linear}
\end{equation}
Thus each $\mathcal C_i$ is a rational polytope. Each $\mathcal C_i$ is nonempty because it contains
$\bm\lambda$ with $\lambda_i=1$ and $\lambda_j=0$ for $j\neq i$.
For every $\bm\lambda\in\Delta_m$, the point $\bm V\bm\lambda$
has a nearest vertex $\bm v_i$, so $\bm\lambda\in\mathcal C_i$.
Thus $\Delta_m=\bigcup_{i=0}^m\mathcal C_i$.

On $\mathcal C_i$, the minimum distance equals
$f_i(\bm\lambda):=\|\bm V\bm\lambda-\bm v_i\|_2^2$.  Hence
\begin{equation}
 R^*(S)=\max_{0\leq i\leq m}\max_{\bm\lambda\in\mathcal C_i}
                           f_i(\bm\lambda).
 \label{eq:svbc-voronoi-decomposition}
\end{equation}
The Hessian $2\bm V^{\mathsf T}\bm V$ is positive semidefinite,
so $f_i$ is convex. 

\textbf{Step 2: rational attainment and exact enumeration.}
A convex function on a compact polytope attains its maximum at a
vertex: if a point is a convex combination of vertices, its value
is at most the corresponding average of their values.
By \cref{eq:svbc-voronoi-decomposition}, some vertex
$\bm\lambda^*$ of a polytope $\mathcal C_i$ yields a global
optimizer $\bm x^*=\bm V\bm\lambda^*$.

Within $\bm1_{m+1}^{\mathsf T}\bm\lambda=1$, each cell has
$2m+1$ inequalities: $m+1$ nonnegativity constraints and $m$
distance comparisons.  At a vertex $\bm\lambda$ of $\mathcal C_i$, one can choose
$m$ active inequalities whose normals, together with
$\bm1_{m+1}$, are linearly independent.
These equalities and $\bm1_{m+1}^{\mathsf T}\bm\lambda=1$
form a nonsingular rational system that uniquely determines
$\bm\lambda$.
Since its coefficients have polynomial encoding length,
Cramer's rule and Hadamard's bound imply that
$\bm\lambda$ and $\bm x=\bm V\bm\lambda$ also have polynomial
encoding length.

We therefore enumerate all subsets of $m$ inequalities in each
cell, impose them as equalities together with
$\bm1_{m+1}^{\mathsf T}\bm\lambda=1$, and discard singular systems
and infeasible solutions.  The active-constraint characterization
ensures that this includes every vertex, including degenerate ones.
Across all cells, the number of systems is at most $ N_m=(m+1)\binom{2m+1}{m}=2^{O(m)}$.
Each system can be solved and checked in polynomial time.
Evaluate $f_i(\bm\lambda)$ at each retained vertex
$\bm\lambda$ of $\mathcal C_i$ and select a pair $(i,\bm\lambda)$
with the largest value.
The point $\bm x=\bm V\bm\lambda$ is then a global optimizer.

\textbf{Step 3: NP membership.}
For a YES-instance, Step~2 provides a rational optimal point
$\bm x^*=\bm V\bm\lambda^*$ with polynomial encoding length.
Use $\bm\lambda^*$ as a certificate.
For a proposed certificate $\bm\lambda$, the verifier checks
\begin{equation}
 \bm\lambda\geq\bm0_{m+1},\qquad
 \bm1_{m+1}^{\mathsf T}\bm\lambda=1,\qquad
 \|\bm V\bm\lambda-\bm v_j\|_2^2\geq\tau
 \quad(0\leq j\leq m)
 \label{eq:svbc-np-verification}
\end{equation}
using exact rational arithmetic in polynomial time.
Thus SVBC-DEC belongs to $\mathrm{NP}$, completing the proof.
\end{proof}

Combining \cref{thm:algorithmic} with the strong NP-hardness for
full-dimensional integer simplices established in
\cref{sec:svbc-hardness} proves \cref{thm:main}.

\section{Conclusion}
\label{sec:conclusion}

We established that $\ell_1$-WMD-DEC is strongly NP-complete and that
SVBC-DEC remains strongly NP-complete for full-dimensional integer
simplices.  Complementary exact algorithms establish polynomial-time
solvability in fixed dimension and produce polynomial-length optimal
solutions: rational for SVBC and quadratic algebraic for weighted
$\ell_1$-ball dispersion.  These representations also provide the
certificates needed to establish NP membership. Whether strong NP-hardness persists for the
unweighted unit $\ell_1$-ball remains an intriguing open question.

\appendix

\section{Auxiliary computations for one-quadric optimization}
\label{app:one-quadric-details}

This appendix supplies the two computational details used in
\cref{thm:one-quadric-optimization}: finding a zero of a concave
quadratic on a rational polytope, and comparing candidate values
exactly.  We keep the notation of the proof so that these procedures
can be applied directly to each enumerated affine slice.

\subsection{Quadratic feasibility on a polytope}
\label{app:quadric-feasibility}

Consider a nonempty rational polytope
\[
 P_J=\{\bm u\in\mathbb R^p:\bm D\bm u\leq\bm\xi\},
 \qquad
 q_J(\bm u)=\gamma_J+2\bm g_J^{\mathsf T}\bm u
                       -\bm u^{\mathsf T}\bm H_J\bm u,
 \qquad \bm H_J\succeq\bm0_p.
\]
We show how to compute rational minimizers $\bm u^{\min}$ and maximizers $\bm u^{\max}$ of $q_J$
on $P_J$, and then use them to find $\bm u\in P_J$ satisfying $q_J(\bm u)=0$, if such a point exists.
The case $p=0$ is a direct evaluation, so assume $p\geq1$.
There are $M$ inequalities, and all input data have polynomial
encoding length in the input length $L$ of
\cref{thm:one-quadric-optimization}.

\paragraph{Computing the minimum}
Concavity reduces this calculation to the vertices $\bm v_j$ of $P_J$.
Indeed, writing any $\bm u\in P_J$ as a convex combination
$\bm u=\sum_j\theta_j\bm v_j$ of vertices gives
\[
 q_J(\bm u)\geq\sum_j\theta_j q_J(\bm v_j)
             \geq\min_j q_J(\bm v_j).
\]
Since the vertices belong to $P_J$, it follows that $h_{\min}=\min_j q_J(\bm v_j) $.

To generate these vertices $\bm v_j$, enumerate all
$I\subseteq\{1,\ldots,M\}$ with
$|I|=\operatorname{rank}\bm D_I=p$, solve
$\bm D_I\bm u=\bm\xi_I$, and retain the solutions satisfying
$\bm D\bm u\leq\bm\xi$.  Every vertex has $p$ independent active
constraint normals, even if $P_J$ is lower dimensional, so this
enumeration includes all vertices.  Evaluating $q_J$ at the retained vertices yields
$h_{\min}$ and a rational minimizer $\bm u^{\min}$.

\paragraph{Computing the maximum}
We compute $h_{\max}$ by enumerating linear systems derived
from the first-order optimality condition.
Since $P_J$ is compact, let $\bm u^\star$ maximize $q_J$ on
$P_J$, and define
\[
 \mathcal A(\bm u^\star)
 =\{i:\bm D_i\bm u^\star=\xi_i\}.
\]
Since $\bm u^\star$ maximizes $q_J$ on the convex set $P_J$,
the first-order optimality condition gives
\[
 \nabla q_J(\bm u^\star)^{\mathsf T}
 (\widehat{\bm u}-\bm u^\star)\leq0
 \qquad\text{for every }\widehat{\bm u}\in P_J.
\]
Thus $\nabla q_J(\bm u^\star)=2\bm g_J-2\bm H_J\bm u^\star$
belongs to the normal cone of $P_J$ at $\bm u^\star$.
Since this cone is generated by the active constraint normals,
\[
 2\bm g_J-2\bm H_J\bm u^\star
 \in\operatorname{cone}
       \{\bm D_i^{\mathsf T}:i\in\mathcal A(\bm u^\star)\}.
\]
By the conic Carath\'eodory theorem, this gradient $\nabla q_J(\bm u^\star)$ is a
nonnegative combination of at most $p$ vectors from
$\{\bm D_i^{\mathsf T}:i\in\mathcal A(\bm u^\star)\}$.
Thus, for some $I\subseteq\mathcal A(\bm u^\star)$ with
$|I|\leq p$, the pair $(\bm u^\star,\bm\mu)$ satisfies
the following linear system in $(\bm u,\bm\mu)$:
\begin{equation}
 \bm D\bm u\leq\bm\xi,\qquad
 \bm D_I\bm u=\bm\xi_I,\qquad
 \bm\mu\geq\bm0,\qquad
 2\bm H_J\bm u-2\bm g_J+\bm D_I^{\mathsf T}\bm\mu=\bm0_p.
 \label{eq:one-quadric-kkt}
\end{equation}
We enumerate all subsets $I$ of size at most $p$, including the
empty set, and test these systems for rational linear feasibility.
The preceding argument guarantees that at least one is feasible.

Conversely, if $(\bm u,\bm\mu)$ satisfies
\cref{eq:one-quadric-kkt}, then $\bm u$ maximizes $q_J$ on $P_J$.
Indeed, for every $\widehat{\bm u}\in P_J$,
$\widehat{\bm u}\in P_J$,
\begin{align}
 q_J(\widehat{\bm u})
 &=q_J(\bm u)
   +\nabla q_J(\bm u)^{\mathsf T}(\widehat{\bm u}-\bm u)
   -(\widehat{\bm u}-\bm u)^{\mathsf T}
      \bm H_J(\widehat{\bm u}-\bm u)\notag\\
 &\leq q_J(\bm u)
   +\bm\mu^{\mathsf T}
      (\bm D_I\widehat{\bm u}-\bm D_I\bm u)
 \leq q_J(\bm u).
 \label{eq:one-quadric-kkt-sufficiency}
\end{align}
A nonempty rational linear system has a rational feasible solution
of polynomial encoding length.  Thus any feasible enumerated system
supplies a rational maximizer $\bm u^{\max}$ and the rational value
$h_{\max}=q_J(\bm u^{\max})$.

\paragraph{Constructing a zero and bounding the cost}
By \cref{eq:one-quadric-interval-test}, no feasible zero exists
if $h_{\min}>0$ or $h_{\max}<0$.
If $h_{\min}=0$, return $\bm u^{\min}$; if $h_{\max}=0$,
return $\bm u^{\max}$.
Otherwise, $h_{\min}<0<h_{\max}$, and we consider
\[
 \bm u(\theta)=(1-\theta)\bm u^{\min}+\theta\bm u^{\max},
 \qquad 0\leq\theta\leq1.
\]
The segment lies in $P_J$, and
$q_J(\bm u(0))=h_{\min}<0<h_{\max}=q_J(\bm u(1))$.
Hence $q_J(\bm u(\theta))=0$ has a solution
$\theta^*\in(0,1)$.  Its degree is at most two and its coefficients
are rational.  The quadratic formula, or a linear equation if the
degree drops, therefore gives an exact root of algebraic degree at
most two.  All coordinates of $\bm u(\theta^*)$ belong to
$\mathbb Q(\theta^*)$.

The vertex and active-set enumerations each have size $M^{O(p)}$.
Every linear solve, feasibility test, and evaluation has polynomial
bit cost.  Cramer's rule and Hadamard's bound give polynomial
encoding length for the rational points \(\bm u^{\min},\bm u^{\max}\) above, and hence for the
coefficients and root representation of \(q_J(\bm u(\theta))=0\).  This proves the required
$M^{O(p)}\operatorname{poly}(L)$ bound and supplies a
polynomial-length point whenever $P_J\cap\{q_J=0\}$ is nonempty.

\subsection{Exact comparisons of quadratic algebraic numbers}
\label{app:quadratic-arithmetic}

The candidate values may be irrational, but each has the form
$a+b\sqrt{\Delta}$ with $a,b,\Delta\in\mathbb Q$ and $\Delta\geq0$.
Their comparison reduces to rational comparisons by checking
signs before squaring.

For $s=e+f\sqrt{\Delta}$ with rational $e,f,\Delta$ and
$\Delta\geq0$, handle zero terms first.  If the two nonzero terms
have the same sign, that is the sign of $s$.  If their signs differ,
compare $e^2$ with $f^2\Delta$.  Equality gives $s=0$; otherwise
$s$ has the sign of the term with larger absolute value.
This determines the sign exactly with one squaring comparison.

Now consider
\[
 \alpha_1=a_1+b_1\sqrt{\Delta_1},
 \qquad
 \alpha_2=a_2+b_2\sqrt{\Delta_2},
 \qquad \Delta_1,\Delta_2\geq0,
\]
with rational coefficients.  Comparing them is equivalent to
comparing
\[
 x=(a_1-a_2)+b_1\sqrt{\Delta_1},
 \qquad y=b_2\sqrt{\Delta_2}.
\]
The preceding test determines the sign of $x$.
The sign of $y$ is zero if $b_2=0$ or $\Delta_2=0$, and otherwise
is the sign of $b_2$.  If either number is zero or their signs differ,
their order follows immediately.  If both are positive, compare
$x^2$ and $y^2$; if both are negative, reverse that comparison.
The difference of the squares is
\[
 x^2-y^2
  =\bigl((a_1-a_2)^2+b_1^2\Delta_1-b_2^2\Delta_2\bigr)
     +2(a_1-a_2)b_1\sqrt{\Delta_1},
\]
so its sign is determined by the single-radical test.

Thus checking signs before squaring reduces the comparison of
any pair of candidate values to comparisons of rational numbers.  The fixed number of rational operations
preserves polynomial encoding length, and each comparison takes
polynomial time.

\section{Rational realization of the balanced spherical cluster}
\label{app:svbc-realization}

We prove \cref{lem:svbc-realization} by constructing the prescribed
norm vectors $\bm g,\bm h,\bm h_0$, a basis of $W_0$ consisting
of integer vectors, and the balanced cluster $\bm c_0,\ldots,\bm c_r$.  The final estimates bound the common
denominator numerically, as required for strong NP-hardness.

\begin{proof}[Proof of \cref{lem:svbc-realization}]
\textbf{Step 1: realizing the prescribed squared norms.}
We first construct the tag vector $\bm g$ and the height vectors
$\bm h,\bm h_0$ with denominator $q$.
By the four-square theorem, choose
$\widehat{\bm g},\widehat{\bm h},\widehat{\bm h}_0\in\mathbb Z^4$
such that
\begin{equation}
 \|\widehat{\bm g}\|_2^2=q(q-2),\qquad
 \|\widehat{\bm h}\|_2^2=q(4q-5),\qquad
 \|\widehat{\bm h}_0\|_2^2=4q^2-5q-1.
 \label{eq:svbc-four-square-tags}
\end{equation}
Then $ \bm g={\widehat{\bm g}}/q$, $ \bm h={\widehat{\bm h}}/q$ and $\bm h_0={\widehat{\bm h}_0}/q $
have precisely the squared norms required in
\cref{eq:svbc-parameters,eq:svbc-lowered-anchor}.

These norm identities imply that the coordinates of
$\bm g,\bm h,\bm h_0$ have magnitude $O(1)$.
We next construct a basis of $W_0$ consisting of integer vectors.

\textbf{Step 2: a basis of $W_0$ consisting of integer vectors.}
We construct a basis of $W_0$ with polynomially bounded integer
entries.  Its explicit size bound will control the common
denominator in the next step.  By definition,
\[
 W_0=\mathcal L^\perp
 =\{\bm z\in\mathbb R^{d_0}:
 \bm z^{\mathsf T}(\bm p_i-\bar{\bm a})=
 \bm z^{\mathsf T}(\bm w_i-\bar{\bm a})=0,
  1\leq i\leq n\}.
\]
In the identity, incidence, tag, and height blocks, write
\[
 \bm z=\bm x\oplus\bm y\oplus
 (\bm\tau_1\oplus\cdots\oplus\bm\tau_n)\oplus\bm z_h,
\]
where $\bm x\in\mathbb R^n$, $\bm y\in\mathbb R^{3q}$,
$\bm\tau_i\in\mathbb R^4$, and $\bm z_h\in\mathbb R^4$.

For $1\leq i\leq n$, blockwise inner products give
\begin{align*}
 \bm z^{\mathsf T}(\bm p_i-\bar{\bm a})
 =
 x_i+(\bm b_i-\bm u)^{\mathsf T}\bm y
 -\bm h_0^{\mathsf T}\bm z_h,\quad
 \bm z^{\mathsf T}(\bm w_i-\bar{\bm a})=
 (\bm b_i-\bm u)^{\mathsf T}\bm y
 +\bm g^{\mathsf T}\bm\tau_i
 -\bm h_0^{\mathsf T}\bm z_h.
\end{align*}
Multiplying these equalities by $q$ and using
$q\bm u=\bm1_{3q}$,
$q\bm g=\widehat{\bm g}$, and
$q\bm h_0=\widehat{\bm h}_0$, we see that
$\bm z\in W_0$ is equivalent to
\begin{align}
qx_i+(q\bm b_i-\bm1_{3q})^{\mathsf T}\bm y
-\widehat{\bm h}_0^{\mathsf T}\bm z_h&=0,
\label{eq:svbc-complement-p}\\
(q\bm b_i-\bm1_{3q})^{\mathsf T}\bm y
+\widehat{\bm g}^{\mathsf T}\bm\tau_i
-\widehat{\bm h}_0^{\mathsf T}\bm z_h&=0.
\label{eq:svbc-complement-w}
\end{align}

To construct the tag directions, set
$\bm y=\bm0_{3q}$ and $\bm z_h=\bm0_4$. Then \cref{eq:svbc-complement-p} and \cref{eq:svbc-complement-w} require
$\bm x=\bm0_n$ and
$\widehat{\bm g}^{\mathsf T}\bm\tau_i=0$
for every $i$. Since $\|\widehat{\bm g}\|_2^2=q(q-2)>0$,
$\widehat{\bm g}$ is nonzero.  Without loss of generality, assume
that $\widehat g_1\neq0$.  For $\ell=2,3,4$, set
\begin{equation}
 \bm\kappa_\ell
 =\widehat g_1\bm e_\ell^{(4)}
  -\widehat g_\ell\bm e_1^{(4)}\in\mathbb Z^4,
 \qquad \text{so}\quad
 \widehat{\bm g}^{\mathsf T}\bm\kappa_\ell
 =\widehat g_1\widehat g_\ell-\widehat g_\ell\widehat g_1=0.
 \label{eq:svbc-tag-kernel}
\end{equation}
Ones can check that
$\bm\kappa_2,\bm\kappa_3,\bm\kappa_4$ are linearly independent. For $1\leq i\leq n$ and $\ell\in\{2,3,4\}$, define
$\bm v_{i,\ell}^{\mathrm{tag}}\in\mathbb Z^{d_0}$ by
\begin{equation}
 \bm v_{i,\ell}^{\mathrm{tag}}
 =\bm0_n\oplus\bm0_{3q}\oplus
  (\bm\tau_1\oplus\cdots\oplus\bm\tau_n)\oplus\bm0_4,
 \qquad
 \bm\tau_j=
 \begin{cases}
  \bm\kappa_\ell,&j=i,\\
  \bm0_4,&j\neq i.
 \end{cases}
 \label{eq:svbc-tag-basis}
\end{equation}
Equations
\cref{eq:svbc-complement-p,eq:svbc-complement-w,eq:svbc-tag-kernel}
show that $\bm v_{i,\ell}^{\mathrm{tag}}\in W_0$.
Vectors $\bm v_{i,\ell}^{\mathrm{tag}}$ associated with different values of $i$ have disjoint
tag-block supports.  Hence \cref{eq:svbc-tag-basis} gives $3n$
linearly independent integer vectors in $W_0$.
These vectors account for the $3n$ degrees of freedom in the tag
blocks orthogonal to $\widehat{\bm g}$.  Since $ \dim W_0=3n+3q+4 $,
it remains to construct $3q$ independent directions associated with
the incidence block $\bm y$ and four associated with the height block
$\bm z_h$. 

To construct the incidence directions, set
$N_g=\|\widehat{\bm g}\|_2^2$.
For each $j\in\{1,\ldots,3q\}$, define
$\bm v_j^{\mathrm{inc}}\in\mathbb Z^{d_0}$
by the following coordinate blocks:
\begin{equation}
\begin{aligned}
 \bm y&=qN_g\bm e_j,       &\qquad \bm z_h&=\bm0_4,\\
 x_i&=-(qB_{ji}-1)N_g,     &
 \bm\tau_i&=-q(qB_{ji}-1)\widehat{\bm g}
 &&(1\leq i\leq n).
\end{aligned}
\label{eq:svbc-incidence-basis}
\end{equation}

For each height coordinate $\ell\in\{1,2,3,4\}$, define
$\bm v_\ell^{\mathrm{ht}}\in\mathbb Z^{d_0}$ by
\begin{equation}
\begin{aligned}
 \bm z_h&=qN_g\bm e_\ell,  &\qquad \bm y&=\bm0_{3q},\\
 x_i&=N_g\widehat h_{0,\ell},&
 \bm\tau_i&=q\widehat h_{0,\ell}\widehat{\bm g}
 &&(1\leq i\leq n).
\end{aligned}
\label{eq:svbc-height-basis}
\end{equation}

Substitution into
\cref{eq:svbc-complement-p,eq:svbc-complement-w}, using
$N_g=\|\widehat{\bm g}\|_2^2$, shows that
\[
 \bm v_j^{\mathrm{inc}},\ \bm v_\ell^{\mathrm{ht}}
 \in W_0\cap\mathbb Z^{d_0}
 \qquad
 (1\leq j\leq3q,\ 1\leq\ell\leq4).
\]

Altogether, we have constructed the family $\mathcal V_0\subseteq W_0\cap\mathbb Z^{d_0}$ as
\begin{equation}
\mathcal V_0
=
\left\{
 \bm v_{i,\ell}^{\mathrm{tag}}:
 1\leq i\leq n,\ \ell=2,3,4
\right\}
\cup
\left\{
 \bm v_j^{\mathrm{inc}}:
 1\leq j\leq3q
\right\}
\cup
\left\{
 \bm v_\ell^{\mathrm{ht}}:
 1\leq\ell\leq4
\right\}.
\label{eq:svbc-W0-basis-family}
\end{equation}

The family $\mathcal V_0$ contains
$s_0=3n+3q+4=\dim W_0$ vectors.
In any vanishing linear combination, projection onto the incidence
and height blocks forces all coefficients of
$\bm v_j^{\mathrm{inc}}$ and $\bm v_\ell^{\mathrm{ht}}$ to vanish.
The remaining coefficients vanish because the tag vectors are
linearly independent.
Thus $\mathcal V_0$ is a basis of $W_0$;
enumerate its vectors as $\bm v_1,\ldots,\bm v_{s_0}$. 

We next bound the basis vectors to control $\Gamma$.
The defining norm identities imply
\[
 \|\widehat{\bm g}\|_\infty=O(q),\qquad
 N_g=q(q-2)=O(q^2),\qquad
 \|\widehat{\bm h}_0\|_\infty=O(q).
\]

Substituting these bounds into
\cref{eq:svbc-tag-basis,eq:svbc-incidence-basis,eq:svbc-height-basis}
shows that every basis vector has coordinates of magnitude $O(q^3)$.
Each vector has at most $O(n+q)$ nonzero coordinates.  Hence
\begin{equation}
 \max_{1\leq i\leq s_0}\|\bm v_i\|_\infty=O(q^3),
 \qquad
 \max_{1\leq i\leq s_0}\|\bm v_i\|_1
 =O((n+q)q^3).
 \label{eq:svbc-basis-bound}
\end{equation}

Recall that $s_0=\dim W_0$ and
$r=\dim W=s_0+9$, since $W=W_0\oplus\mathbb R^9$. We now use the basis $\bm v_1,\ldots,\bm v_{s_0}$ of $W_0$
to construct the balanced equal-norm cluster
$\bm c_0,\ldots,\bm c_r\in W$.

\textbf{Step 3: a balanced equal-norm cluster.}
Split the nine additional coordinates as
$\mathbb R^9=\mathbb R^4\oplus\mathbb R^5$.
Our goal is to construct $\bm c_0,\ldots,\bm c_r$
satisfying \cref{eq:svbc-cluster-properties}.

Let $\bm v_\Sigma=\sum_{i=1}^{s_0}\bm v_i$ and choose
\begin{equation}
 \Gamma
 =3+\|\bm v_\Sigma\|_1+\max_i\|\bm v_i\|_1.
 \label{eq:svbc-Gamma}
\end{equation}
By the four-square theorem, choose nonnegative integer vectors
$\bm t_i,\bm t_0\in\mathbb Z^4$ satisfying
\begin{align}
 \|\bm t_i\|_2^2
 =\Gamma^2-\|\bm v_i\|_2^2
 \quad(1\leq i\leq s_0),\quad
 \|\bm t_0\|_2^2 =\Gamma^2-\|\bm v_\Sigma\|_2^2-5.
 \label{eq:svbc-four-square-cluster}
\end{align}
Because every right-hand side is positive and $s_0\geq4$, permute the
coordinates of $\bm t_1,\ldots,\bm t_4$ so that the $k$th coordinate
of $\bm t_k$ is positive.

In
$W=W_0\oplus\mathbb R^4\oplus\mathbb R^5$, define
\begin{align}
 &\bm u_i=(\bm v_i,\bm t_i,\bm0_5)\quad(1\leq i\leq s_0),\quad
 \bm u_0=(-\bm v_\Sigma,\bm t_0,\bm1_5), \notag\\
 &\bm d_j=(\bm0_{d_0},-\Gamma\bm e_j^{(9)})
          \quad (1\leq j\leq9).
 \label{eq:svbc-integer-cluster}
\end{align}
Every vector in \cref{eq:svbc-integer-cluster} has norm $\Gamma$. $\bm{e}_j^{(9)}$ is the $j$-th standard basis vector in $\mathbb{R}^9$. 

Moreover, define
\begin{equation}
 \bm T
 =\left(\bm t_0+\sum_{i=1}^{s_0}\bm t_i,\bm1_5\right)
 \in\mathbb Z_{+}^9.
\end{equation}
The corresponding vector
$\bm t_k$ has a positive $k$th entry and every $\bm t_i$ is
componentwise nonnegative.
Hence $\bm T$ is strictly positive.  It follows directly from
\cref{eq:svbc-integer-cluster} that
\begin{equation}
 \bm u_0+\sum_{i=1}^{s_0}\bm u_i+
 \sum_{j=1}^9\frac{T_j}{\Gamma}\bm d_j=\bm0_d.
 \label{eq:svbc-positive-cluster-relation}
\end{equation}
All coefficients in this linear relation are strictly positive.

By construction, the $r=s_0+9$ vectors
$\bm u_1,\ldots,\bm u_{s_0},\bm d_1,\ldots,\bm d_9$
form a basis of $W$, as seen by projection onto $W_0$
and the final $\mathbb R^9$ block.
Together with the strictly positive relation
\cref{eq:svbc-positive-cluster-relation}, this implies that
$\bm u_0,\bm u_1,\ldots,\bm u_{s_0},\bm d_1,\ldots,\bm d_9$
are affinely independent and that $\bm0_d$ lies in the relative
interior of their convex hull.

Define
\[
\bm c_0=\frac{\bm u_0}{q\Gamma},
\qquad
\bm c_i=\frac{\bm u_i}{q\Gamma}
\quad(1\leq i\leq s_0),
\qquad
\bm c_{s_0+j}=\frac{\bm d_j}{q\Gamma}
\quad(1\leq j\leq9).
\]
The common scaling preserves affine
independence and the relative-interior property, while
\[
\|\bm c_j\|_2=\frac1q=\varepsilon,
\qquad
q\Gamma\bm c_j\in\mathbb Z^d
\quad(0\leq j\leq r).
\]
Hence \cref{eq:svbc-cluster-properties} holds.

\textbf{Step 4: construction time and a common denominator.}
We finish by bounding $\Gamma$ and showing that the four-square
representations can be found in polynomial time.
Since $s_0=O(n+q)$, \cref{eq:svbc-basis-bound} gives
\[
 \|\bm v_\Sigma\|_1
 \leq s_0\max_i\|\bm v_i\|_1
 =O((n+q)^2q^3).
\]
The definition \cref{eq:svbc-Gamma} therefore gives $ q\Gamma\bm c_j\in\mathbb Z^d$ and $\Gamma=O((n+q)^2q^3) $.

For each integer $N$ in
\cref{eq:svbc-four-square-tags,eq:svbc-four-square-cluster},
enumerate $a,b,c\in\{0,\ldots,\lfloor\sqrt N\rfloor\}$
and test whether $N-a^2-b^2-c^2$ is a nonnegative square.
Lagrange's four-square theorem guarantees success.
Since each $N$ is bounded by $O(q^2)$ or by $\Gamma^2$,
the $O(N^{3/2})$ tests take polynomial time in the explicit
X3C input size.

Thus $\bm g,\bm h,\bm h_0$ and $\bm c_0,\ldots,\bm c_r$
can be constructed in polynomial time with the properties
stated in \cref{lem:svbc-realization}.
\end{proof}

\end{document}